\documentclass[11pt]{amsart}

\usepackage{mathrsfs}
\usepackage{amsmath,amsthm,amsfonts,amssymb,mathrsfs,amscd,mathtools}
\usepackage{graphicx}
\usepackage{dsfont}
\usepackage{hyperref}
\usepackage[svgnames]{xcolor}
\usepackage{pgffor} 
\usepackage{pgfplots}
\usepackage{tikz-cd}
\usepackage[mathcal]{euscript}
\usepackage{subcaption}
\usepackage[english]{babel}
\usepackage{url}
\hypersetup{
    colorlinks=true,
    linkcolor=Blue, 
    citecolor=magenta,
}
\usepackage{extarrows}
\newtheorem{maintheorem}{Theorem}

\newcommand{\C}{\mathbb{C}}
\newcommand{\R}{\mathbb{R}}
\newcommand{\N}{\mathbb{N}}
\newcommand{\Z}{\mathbb{Z}}
\newcommand{\Q}{\mathbb{Q}}
\renewcommand{\P}{\mathbb{P}}
\newcommand{\ohne}{\smallsetminus}
\renewcommand{\O}{\mathcal{O}}
\renewcommand{\P}{\mathbb{P}}
\DeclareMathOperator{\Spec}{Spec}
\DeclareMathOperator{\pr}{pr}
\DeclareMathOperator{\id}{id}
\DeclareMathOperator{\CH}{CH}
\DeclareMathOperator{\Pic}{Pic}
\newcommand{\X}{\mathcal{X}}		
\newcommand{\num}{\mathrm{num}}
\newcommand{\alg}{\mathrm{alg}}
\newcommand{\Zykel}{\mathrm{Z}}
\newcommand{\PP}{\mathcal{P}}
\newcommand{\NN}{\mathcal{N}}
\newcommand{\AJ}{\mathrm{AJ}}
\newcommand{\E}{\mathcal{E}}		
\newcommand{\A}{\mathscr{A}}

\newcommand{\an}{\mathrm{an}}
\DeclareMathOperator{\PD}{PD}
\newcommand{\BB}{\mathcal{B}}
\newcommand{\Sbar}{\overline{S}}
\newcommand{\Zykle}{\mathrm{Z}}
\newcommand{\SCH}{\mathcal{CH}}
\newcommand{\calA}{\mathcal{A}}
\newcommand{\calB}{\mathcal{B}}
\newcommand{\bG}{\mathbb{G}}
\newcommand{\bE}{\mathbb{E}}
\renewcommand{\hom}{\mathrm{hom}}

\calclayout

\newtheorem{thm}{Theorem}[section]
\newtheorem{cor}[thm]{Corollary}
\newtheorem{prop}[thm]{Proposition}
\newtheorem{lem}[thm]{Lemma}
\newtheorem{conj}[thm]{Conjecture}
\newtheorem{quest}[thm]{Question}

\theoremstyle{definition}
\newtheorem{defn}[thm]{Definition}

\newtheorem{ex}[thm]{Example}

\newtheorem{rem}[thm]{Remark}
\newtheorem{void}[thm]{}

\theoremstyle{remark}
\newtheorem{conv}[thm]{Convention}
\newtheorem{notn}[thm]{Notation}

\numberwithin{equation}{subsection}
\pgfplotsset{compat=1.9}
\title[]{Degeneration of the archimedean height pairing of algebraically trivial cycles}
\author[]{Zhelun Chen}
\date{\today}
\thanks{2020 Mathematical Subject Classification: 14G40, 14C30, 14C25}
\address{Mathematical Institute, Leiden University, Einsteinweg 55, 2333 CC Leiden, The Netherlands}
\email{zl.chen1729@gmail.com}
\begin{document}

\begin{abstract}
We consider  the limiting behaviour of the  archimedean height pairing for homologically trivial  algebraic cycles  in a degenerating one-parameter family of smooth projective complex varieties. We conjecture that the limit is controlled by the non-archimedean geometric height pairing of the cycles on the generic fiber and verify this for algebraically trivial cycles, assuming a  conjecture of Griffiths on incidence equivalence.  Our work offers a more geometric understanding of a related asymptotic result of Brosnan--Pearlstein and suggests a new perspective on the positivity of the Beilinson--Bloch height pairing over a one-variable complex function field.
\end{abstract}

\maketitle

\tableofcontents
\section{Introduction}
The principal motivation of the present manuscript is to understand an asymptotic result on the archimedean height pairing due to Brosnan--Pearlstein \cite{BP19} in the geometric setting. We obtain a satisfactory answer for algebraically trivial cycles, conditional on a conjecture of Griffiths on incidence equivalence of algebraic cycles   (Conjecture \ref{Griffiths conjecutre for incidentally trivial cycles}).

Let us introduce the problem properly and outline the main ingredients of our solution. 
\subsection{The problem}
\begin{void}[Standard setup]\label{standard setup}
		Let $f\colon X\to S$ be a smooth projective morphism of complex algebraic varieties of fiber dimension $d$ over a smooth complex algebraic curve $S$; let $\eta\in S$ be the generic point. Let $p, q\in\N$ such that $p+q=d+1$. Let $Z\in \Zykle^p(X/S),\, W\in \Zykle^q(X/S)$ be relative (flat) algebraic cycles of relative codimension $p, q$ respectively. Suppose the cycles $Z, W$ do not intersect over~$S$, i.e.\ their supports on $X$ do not intersect: $|Z|\cap |W|=\emptyset$. Moreover, assume that the fibers $Z_s$ and $W_s$ ($s\in S(\C)$) are  homologically trivial. 
\end{void} 

In this setup, the fiberwise \textbf{archimedean  height pairing} of $Z, W$ is defined to be
	\begin{equation}
		S(\C)\to \R, \quad s\mapsto h_\infty(s)\coloneq \langle Z_s, W_s\rangle_\infty\coloneq -\int_{Z_s(\C)}g_{W}\in \R,
	\end{equation}  where $g_{W}$ denotes a suitably normalized  Green's current associated to $W$ which is a smooth form away from the support $|W|$. The assumptions in Setup \ref{standard setup} ensure the well-definedness of  the function $h_\infty$  on $S(\C)$, see e.g.\ \cite[\S 1]{Bost}.
\begin{void}
		We are interested in the  limiting behaviour  of $h_\infty$. More precisely, consider a family $f\colon  X\to  S$  as in \ref{standard setup}. By a \textbf{degeneration} of $X/S$ we mean that there exists a pullback diagram
		\[
		\begin{tikzcd}
			X\ar[hook]{r} \ar[swap]{d}{f} & \overline{X}\ar{d}{\bar{f}}\\
			S\ar[hook]{r} & \overline{S}
		\end{tikzcd}
		\] where $\overline{X}$ is a smooth projective complex variety, $\Sbar$ is a smooth algebraic curve containing~$S$ as a Zariski open subset and $\bar{f}$ is a proper flat morphism; we call $\Sbar$  a (smooth) \textbf{partial compactification} of $S$. The degeneration $\overline{X}/\overline{S}$ is called \textbf{semistable}, if   for every $s_0\in \overline{S}\ohne S$, the fiber $\bar{f}^{-1}(s_0)$  is a simple normal crossing divisor on $\overline{X}$. 
\end{void}
One of the main results of Brosnan--Pearlstein, \cite[Theorem 24]{BP19},  can be  specialized as follows.
\begin{thm}[Brosnan--Pearlstein]\label{BP19, Main}
Let $\overline{X}/\Sbar$ be a semistable degeneration of $X/S$ as in \ref{standard setup}. Let $\Delta\subset \C$ be a (small) analytic neighbourhood around $s_0\in \Sbar\ohne S$ such that $s_0$ is mapped to $0\in \Delta$; let $t$  be a holomorphic coordinate of $\Delta$. There is a  \emph{rational} number $\mu_0$  such that the smooth function
\[
h_\infty(t)+\mu_0\log|t|
\] on $\Delta^*$ extends to a continuous function on $\Delta$.
\end{thm}
\begin{rem}
Let us immediately point out that Brosnan--Pearlstein's result holds more generally for  admissible  variations (w.r.t.\ a partial compactification $S\subset \Sbar$) of mixed Hodge structure of biextension type over $S$, which do not necessarily arise from actual algebraic cycles, and the $\mu_0$ reflects the monodromy of the variation about $s_0$. Our Setup \ref{standard setup} gives rise to a such admissible variation by M. Saito's work (\cite{Saito}). According to Hain's philosophy (\cite{Hain90}),  the function $h_\infty$ can be viewed as  the period map of a real variation of mixed Hodge structure, and the admissibility ensures its ``good'' asymptotic behaviour, as confirmed in full generality by Brosnan--Pearlstein.
\end{rem}
The perspective from which we understand the Brosnan–Pearlstein result is crystallized in the following
\begin{conj}\label{Conjecture: geometric height pairing}
		The number $\mu_0$ from Theorem \ref{BP19, Main}  equals  the local geometric height pairing $\langle Z, W\rangle_0$   of the cycles $Z,\, W$ over  $s_0\in \Sbar\ohne S$. 
\end{conj}
\begin{void}[Geometric height pairing]
	One way to realize the (local) geometric height pairing, that is, the (local) Beilinson--Bloch height pairing over a function field, is as follows. First, we have to \emph{assume} that $Z$ and $W$ admit  so-called admissible liftings $Z^h, W^h$ on $\overline{X}$ (\emph{cf}.\ \ref{admissible cycle class}); these are certain $\Q$-cycles which are fiberwise numerically trivial, i.e.\ the intersection  with vertical $\Q$-cycles  on  $\overline{X}/\Sbar$ of complementary dimension is trivial. By assumption, the algebraic intersection product $Z^h\cdot W^h$ is a 0-cycle on~$\overline{X}$ supported on the (possibly) singular fibers of $\bar{f}$.  We then define the \textbf{geometric height pairing} over~$s_0$ as the  multiplicity   
	\[
	\langle Z,W\rangle_0\coloneq \mathrm{mult}_{s_0}(\bar{f}_*(Z^h\cdot W^h))\in\Q
	\] of the 0-cycle  $\bar{f}_*(Z^h\cdot W^h)$ (which is a divisor on $\Sbar$) at $s_0$.  The definition only depends on the generic data and is insensitive to the choices involved.  See   \S \ref{The geometric height pairing à la Kahn} for more details
\end{void}
\begin{rem}[Silverman's limit theorem]
	Conjecture \ref{Conjecture: geometric height pairing} could be considered as a local variant of a classical result of Silverman concerning the limit behaviour of the global Néron--Tate  height pairing, a question proposed by M. Artin, \emph{cf}.\ \cite{Tate83}.  We reformulate \cite[Theorem B]{Sil83} as follows.
	
	Let $A/S$ be an abelian scheme over a smooth curve $S$ defined over a number field $k$; let $\eta\in S$ be the generic point. Then for  every sections $P, Q\in A(S)$ we have
	\[
	\lim\limits_{s\in S(\overline{k}),\, h_{\overline{S}}(s)\to \infty}\frac{\langle P_s, Q_s\rangle_{A_s}}{h_{\overline{S}}(s)}=\langle P_\eta, Q_\eta\rangle_{A_\eta},
	\] where $\langle -,-\rangle_{A_{(-)}}$ denotes the global Néron--Tate height pairing at the corresponding fibers, and $h_{\overline{S}}$ is a Weil height on the smooth projective compactification $\overline{S}$ of $S$.
\end{rem}
\begin{rem}
	We record another related asymptotic result due to B. Harris and B. Wang without assuming the homological triviality, see \cite{HarrisWang}. Let  $Z$, $W$ be closed complex submanifolds of a compact Kähler manifold $X$ of complementary dimension as above. Assume $Z, W$ intersect ``cleanly'', i.e.\ the intersection $M\coloneq Z\cap W$ is a complex submanifold and the tangent bundles satisfy $TZ|_M\cap TW|_M=TM$. Take an analytic deformation $(Z_t)_{t\in \Delta}$ of $Z$ such that $Z_0=Z$. The statement is that the smooth function\footnote{In the reference the asymptotics reads $e\log|t|^2$. We ignore the factor 2 by using a different normalization of the archimedean height pairing.}
	\[
	\langle Z_t,W\rangle_{\infty}+e\log|t|
	\] extends to a continuous function on $\Delta$, where
	\[
	e= \int_M c_m\left(\frac{TX|_M}{TZ|_M+TW|_M}\right)
	\] stands for the $m$-th Euler number of the quotient  tangent bundle (which is still a vector bundle by assumption),  $m=\dim M$ being the complex dimension of $M$. 
\end{rem}
Let us return to  Setup \ref{standard setup} and resume the notation there.

We can reformulate our conjecture in terms of an isomorphism of extensions of certain line bundles.  One  arises naturally from Hodge Theory.
\begin{void}[Biextension bundle]
 Recall   the cycles $Z, W$ are assumed to be homologically trivial, thus we can consider the associated (holomorphic) normal function
\[
(\nu_Z,\nu_W^\vee)\colon S\to J^p_{X/S}\times_S J^{p,\vee}_{X/S},
\] where $J^p_{X/S}$ denotes the relative $p$-th Griffiths intermediate Jacobian, and we identify~$J^q_{X/S}$ with the dual torus $J^{p,\vee}_{X/S}$ using the Poincaré duality. There is a canonical  Poincaré biextension (holomorphic) line bundle $\PP$  over $J^p_{X/S}\times_S J^{p,\vee}_{X/S}$, rigidified along the identity section. The \textbf{biextension bundle} associated to $Z, W$ is given by the pullback
\[
\BB_{Z,W}\coloneq (\nu_Z,\nu_W^\vee)^*\PP\in \Pic(S).
\]
\end{void}
\begin{void}[Biextension height]
	The Poincaré bundle $\PP$ carries a smooth canonical (biextension) metric $||.||$ that encodes the archimedean height. Indeed, the assumption that $Z, W$ do not intersect over $S$  induces a nowhere-vanishing (holomorphic) section $1_{Z,W}\in\Gamma(S,\BB_{Z,W})$ such that
	\begin{equation}
		-\log||1_{Z,W}(t)||=\langle Z_t,W_t\rangle_\infty,\quad t\in S(\C).
	\end{equation} This also shows that the function $h_\infty(t)=\langle Z_t,W_t\rangle_\infty$ is smooth. We denote again by $||.||$ the  metric on $\BB_{Z,W}$ induced by the pullback. See \cite[\S  3]{Hain90} for a more systematic discussion. 
\end{void} 
The following notion is originally due to Mumford (\cite{Mum77}) and used (independently) in D. Lear's thesis \cite{Lear}. 
\begin{defn}
Let $X$ be a complex manifold and $X\subset \overline{X}$ a smooth partial compactification such that the boundary $D\coloneq \overline{X}\ohne X$ is a simple normal crossing divisor. Let $(L,||.||)$ be a smooth metrized holomorphic line bundle over $X$. A \textbf{Lear extension} $[L,||.||]_{\overline{X}}$ of $L$ is a $\Q$-line bundle over $\overline{X}$ such that
\begin{itemize}
	\item  it extends $L$ as a holomorphic line bundle;
	\item the metric $||.||$  extends continuously over $\overline{X}\ohne D^{\mathrm{sing}}$.
\end{itemize}
Here $D^{\mathrm{sing}}\subset \overline{X}$ denotes the singular locus of $D$ which is a closed subset of $\overline{X}$ of codimension at least two.
\end{defn}
Note that there is at most one Lear extension of $(L,||.||)$. And if $\dim X=1$, the second condition simply means that the metric extends continuously over the whole $\overline{X}$.

We can thus reformulate Brosnan--Pearlstein's result:
\begin{thm}[Brosnan--Pearlstein]\label{BP19, Main, geometric}
	The biextension line bundle $\BB_{Z,W}$ admits a Lear extension over a partial compactification $\Sbar$. In fact, the rational number $\mu_0$ of Theorem \ref{BP19, Main} is exactly the multiplicity
	\[
	\mu_0=\mathrm{mult}_{s_0}(1_{Z,W}),
	\] where $1_{Z,W}$ is seen as a meromorphic section of the Lear extension $[\BB_{Z,W},||.||]_{\Sbar}$.
\end{thm}
\begin{notn}
		For an algebraic cycle $C$ we let $c$ denote its Chow cycle class. Since the Abel--Jacobi map factors through Chow cycle classes, one can define the biextension bundle $\BB_{z,w}$ associated to the  classes of $Z, W$ in a similar way. 
\end{notn}
\begin{void}[Height pairing bundle]
On the other hand, the algebraic intersection product is  well-defined  for Chow cycle classes, so we can introduce the line bundle for the geometric height pairing associated to  $z, w$:
	\[
	\E_{<z,w>}\coloneq \O_{\Sbar}(\bar{f}_*(Z^h\cdot W^h))\in \Pic(\Sbar)_\Q,
	\] where $Z^h, W^h$ are admissible liftings of the actual cycles representing $z, w$. The name is justified by the property that if $\Sbar$ is a  smooth projective curve,
	\[
	\deg\E_{<z,w>}=\langle z,w\rangle_{X_\eta}\in\Q
	\] is the geometric Beilinson--Bloch height pairing. The construction is based on Kahn's framework \cite{Kahn}, see \S \ref{The geometric height pairing à la Kahn} for details.
\end{void}
Now our Conjecture \ref{Conjecture: geometric height pairing} can be reformulated as follows. (On taking a projective compactification, we see that the Lear extension is an algebraic line bundle, which in turn implies that it is  algebraic  over any partial compactification.)
\begin{conj}\label{Conjecture: geometric height=Lear extension}
	Let $\overline{X}/\Sbar$ be a semistable degeneration of $X/S$ as in \ref{standard setup} over a partial compactification $\Sbar\supset S$. \emph{Assume} $Z, W$ admit admissible liftings $Z^h, W^h$ on $\overline{X}/\Sbar$. There is a canonical isomorphism of $\Q$-line bundles over $\Sbar$:
	\[
	[\BB_{Z,W},||.||]_{\Sbar}\cong \E_{<z,w>}.
	\]
\end{conj}
Recall that every algebraically trivial cycle is homologically trivial.  By Theorem \ref{Analog of Künnemann's Lemma 8.1}, every algebraically trivial cycle admits an admissible lifting. For the notion of algebraic triviality, see \S\ref{Algebraically trivial cycles revisit} for a reminder and \cite{ACV19Alg} for a thorough treatment.

The main result of this paper is the following
\begin{maintheorem}[Theorem \ref{Main Result, alg. trivial cycles}]\label{Main Theorem}
Suppose Griffiths's conjecture \ref{Griffiths conjecutre for incidentally trivial cycles} holds for codimension-$p$ algebraic cycles on $X$. Then Conjecture \ref{Conjecture: geometric height=Lear extension} is true, without assuming the semistable degeneration.
\end{maintheorem}

\begin{rem}
	The conjectural isomorphism in Conjecture \ref{Conjecture: geometric height=Lear extension} extends the unconditional isomorphism
	\[
	\BB_{Z,W}\cong \bE_{z,w}\in \Pic(S),
	\] proved by Müller-Stach (\cite[Theorem 1]{MS95}), where the right hand side stands for the fiber of Bloch's biextension 
	\[
	\bE\to \SCH^p_\hom(X/S)\times \SCH^q_\hom(X/S)
	\] at $z, w$, now viewed as global section of the Chow sheaves $\SCH^\bullet_\hom(X/S)$ of homologically trivial cycles over $S$, see \cite{BlochBiext} (and \cite[\S\S 2,3]{dJNormalFunction} for an concise overview of Bloch's results in the complex setting). The latter encodes the algebraic intersection pairing on $X/S$. We will return to this notion (in a slightly different manner) in the final Section \ref{Section: G_m-biextensions}.
\end{rem}
We establish the comparison in Theorem \ref{Main Theorem} by relating the line bundles we want to compare to a  line bundle that encapsulates  the geometric Néron--Tate height pairing. 
\begin{void}[Néron--Tate bundle]
	Following Moret-Bailly (\cite[\S III.3]{Ast129}), to geometric points $x\in A(K), y\in A^\vee(K)$  of an abelian variety $A$ defined over the function field $k$ of $\Sbar$, where $K$ denotes an algebraic closure of $k$, one can associate a $\Q$-line bundle $\BB_{<x,y>}$ on $\Sbar$ such that $\deg\BB_{<x,y>}=\langle x,y\rangle_A$ is the geometric Néron--Tate height pairing, when the base $\Sbar$ is a smooth projective curve, based on the theory of Néron models (this is one reason we restrict to the one-dimensional base case). More details can be found in \S\ref{Geometric NT height pairing}.
\end{void}
	
	\subsection{Künnemann's theorem}\label{Strategies}
	\hfil
\begin{void}
	Let $K$ be an algebraic closure of the function field of the generic fiber of $X/S$, and let $X_K$ denote the geometric generic fiber, $d\coloneq \dim X_K$. There exist certain surjective Abel--Jacobi maps
	\[
	\theta^p \colon \CH^p_\alg(X_K)\to \Pic^p_{X_K}(K)
	\] into the so-called $p$-th higher Picard variety $\Pic^p_{X_K}$ developed by H. Saito (\cite{HSaito}). Here $\CH^p_\alg(X_K)$ denotes Chow group of algebraically trivial cycles of codimension $p$ on $X_K$  ($p\in \{1,\ldots,d\}$). These are certain abelian varieties over $K$ such that $\Pic^1_{X_K}$ is the  Picard variety and $\Pic^d_{X_K}$ is the Albanese variety associated to $X_K$. The kernel of $\theta^p$ are exactly the incidentally trivial cycles on $\CH^p_\alg(X)$. We refer to \S \ref{Regular homomorphisms and higher Picard varieties} for more details.
\end{void}

Inspired by a result of Künnemann (\cite{Kue96}), who showed that in the number field case the arithmetic Beilinson--Bloch height pairing of algebraically trivial can be expressed in terms of the Néron--Tate height pairing of their images in the higher Picard varieties (as numerical identity), we manage to prove
	\begin{thm}[Theorem \ref{Kue96, Thm.8.2}]\label{Künnemann's comparison}
		There is a canonical isomorphism of $\Q$-line bundles over~$\Sbar$:
		\[
		{k^p_X}\E_{<z,w>}=\BB_{<\theta^p(z),\lambda^q_X\theta^q(w)>}.
		\]
	\end{thm}
Here $k^p_X$ is some positive integer  ($k^p_X=1$ if $X_\eta$ is a nice curve), and $\lambda^q_X\colon \Pic^q_{X_K}\to \Pic^{p,\vee}_{X_K}$ denotes a certain isogeny (\emph{cf}.\ \ref{Pic of Poincare p-cycle}), and we  have abused notations by viewing $z, w$ also as algebraically trivial cycle class on the geometric generic fiber $X_K$ so that we can apply the previous Abel--Jacobi map $\theta^\bullet$.

To conclude Theorem \ref{Main Theorem} we shall establish
\begin{thm}[Theorem \ref{comparison biextension vs NT-height line bundles}]\label{height pairing bundle vs biextension bundle}
	Assume Griffiths's Conjecture \ref{Griffiths conjecutre for incidentally trivial cycles} for $\CH^p_\alg(X_K)$. There is a canonical isomorphism of $\Q$-line bundles over $\Sbar$:
	\[
	\BB_{<\theta^p(z),\lambda^q_X\theta^q(w)>}=k^p_X[\BB_{z,w},||.||]_{\Sbar}.
	\]
\end{thm}
\subsection{Higher Picard varieties and algebraic intermediate Jacobians}
\hfil

To prove Theorem \ref{height pairing bundle vs biextension bundle}, we need to compare the ``theories of Jacobians'' (in the sense of Lieberman, see \cite{LiebermanSurvey}) modelled by higher Picard varieties resp.\ by algebraic (or Lieberman) intermediate Jacobian. The latter is recently generalized to the relative setting by Achter--Casalaina-Martin--Vial  so that we can consider such Jacobians $J^\bullet_{a,X_K}$ over (an algebraic closure $K$ of) a complex function field. See \S \ref{Algebraic intermediate Jacobians} and the original sources \cite{ACV19}, \cite{ACV23}.
\begin{void}
	More precisely,  consider the Abel--Jacobi map
	\[
	\AJ^\bullet\colon \CH_\alg^\bullet(X_K)\to J^\bullet_{a,X_K}
	\] restricted to algebraically trivially cycles. There is a morphism
	\begin{equation}
		h^\bullet\colon J^\bullet_{a,X_K}\to \Pic^\bullet_{X_K}
	\end{equation}  of abelian varieties so that we can compare the Poincaré bundles on $\Pic^\bullet_{X_K}$ resp.\ on $J^\bullet_{a,X_K}$. The latter is directly related to the biextension bundle in Conjecture \ref{Conjecture: geometric height=Lear extension}, which governs the asymptotics of the archimedean height of algebraically trivially cycles, and the morphism $h^\bullet$ provides the link to the comparison established in Theorem \ref{Künnemann's comparison}.

However, an ideal comparison  requires $h^\bullet$ to be an isogeny, which is the content of Griffiths's conjecture (see Example \ref{Example: Griffiths conj} for some evidences). We remark that the notion of ``relative higher Picard varieties'' is still missing except for the both extreme codimensions (Remark \ref{relative higher Picard var}).
\end{void}
The last two theorems will be proved in   Section \S \ref{Comparison of the (extended) Poincaré bundles}.  A necessary ingredient is an alternative description of the Lear extension using Néron models, which is based on an extension result of the metrized Poincaré bundle:
\begin{thm}
	Let $A\to S$ be an abelian scheme and $\NN(A)^0\to \Sbar$ be the identity component of the Néron model of $A$ over $\Sbar$. Let $(\PP,||.||)$ be the smooth metrized Poincaré  bundle. Then $(\PP,||.||)$ extends to a continuous metrized line bundle over $\NN(A)^0\times_{{\overline{S}}} \NN(A^\vee)^0$.
\end{thm}
The proof uses the formalism of (analytic) adelic line bundles by Yuan--Zhang \cite{YZ21}, see  the comment on Theorem \ref{continuous extension of the Poincaré bundle}. By uniqueness of Lear extensions, we immediately obtain
\begin{cor}[Corollary \ref{Lear extension and Neron model}]\label{Lear extension via Neron model}
	There is a canonical isomorphism of $\Q$-line bundles over~$\Sbar$:
	\[
	[\BB_{z,w},||.||]_{\Sbar}\cong \BB_{<\nu_z,\nu_w^\vee>}.
	\] 
\end{cor}
Here we abusively use $\nu_z$ to denote the $K$-point of the abelian variety  $J^p_{a,X_K}$ induced by the normal function $\nu_z\colon S\to J^p_{a,X/S}$; similarly for $\nu_w$.

To finalize the comparison of the line bundles we need to compare the  morphisms
\[
J^q_{a,X_K}\to J^{p,\vee}_{a,X_K}
\]  induced by the Poincaré duality $J^q_{X/S}\cong J^p_{X/S}$ and   $h^{p,\vee}\circ \lambda^q\circ h^q$ resulted from the comparison of the Jacobians. This amounts to comparing two $\bG_m$-biextensions, a concept due to Mumford (\cite{MumfordBiext}) and Grothendieck, see \S \ref{Section: G_m-biextensions}.

\subsection{Connection to the positivity of the height pairing}\label{Section: positivity of height pairing}
\hfil

Combined with a positivity result of Brosnan--Pearlstein (see below), Conjecture \ref{Conjecture: geometric height=Lear extension} offers a new way to study the positivity of the geometric height pairing. The discussion below will not be pursued further in the main body of the text; nonetheless, it provides an additional motivation for the present work.

Recall that  the Néron--Tate height pairing $\langle,\rangle_A$ on an abelian variety $A$ over a global field $k$ is positive definite. As a generalization of this pairing, the (conjecturally defined) Beilinson--Bloch height pairing $\langle,\rangle_X$  on an arbitrary smooth projective variety $X$ of dimension $d$ over~$k$  is expected to satisfy analogous positivity (see \cite[Conjectures 5.3 \& 5.4]{BeilinsonHP}):
\begin{conj}[Beilinson]\label{Beilinson's conjecture on height pairing}
	Let $L\colon \CH^\bullet(X)\to \CH^{\bullet+1}(X)$ be the Lefschetz operator induced by an ample line bundle $L$ on $X$. 
	\begin{enumerate}
		\item  For $p\leq \frac{d+1}{2}$, the iterated homomorphism 
		\[
		L^{d+1-2p}\colon \CH^p_{\hom}(X)_\Q\to \CH^{d+1-p}_{\hom}(X)_\Q
		\] is an isomorphism. 
		\item Let $z\in \CH^p_{\hom}(X)_\Q$ be a primitive cycle, i.e. $z\in \CH^p_{\hom}(X)_\Q\ohne\{0\}$ and $L^{d+2-2p}z=0$, then 
		\[
		(-1)^p\langle z,L^{d+1-2p}z\rangle_X>0.
		\]
	\end{enumerate}
\end{conj}
\begin{rem}[Dual normal function]\label{dual normal function}
For $z\in \CH^p_{\hom}(X)_\Q$ let us write $$z^\vee\coloneq (-1)^pL^{d+1-2p}z.$$  Beilinson's Conjecture \ref{Beilinson's conjecture on height pairing}(1) says every $w\in \CH^{d+1-p}_{\hom}(X)_\Q$ is of the form $z^\vee$ for some $z\in \CH^p_{\hom}(X)_\Q$.	 Unfolding the definition of Poincaré duality, one can show 
	\begin{equation}
		\nu_z^\vee=\nu_{z^\vee}.
	\end{equation}
\end{rem}
\begin{void}
	Let $X/S$ be a smooth projective family as in Setup \ref{standard setup}; let  $\overline{X}/\Sbar$ be a semistable degeneration with $\Sbar$ being the smooth \emph{projective} compactification of $S$. Let $z\in \CH^p_{\hom}(X/S)$. We also view $z\in \CH^p_{\hom}(X_K)$ as a homologically trivial cycle on the geometric generic fiber $X_K$ of $X/S$.  Suppose   Conjecture \ref{Conjecture: geometric height=Lear extension} holds for this $z$, then
	\[
	L\coloneq \E_{<z,z^\vee>}\cong [\BB_{z,z},||.||]_{\Sbar}.
\] We know from \cite[Theorem 27]{BP19}  that the smooth semipositive biextension metric $h$ on $\BB_{z,z}$ extends to a  continuous semipositive metric on $L$ over $\Sbar$. In other words, 
\[
\deg(L,{h})\coloneq \int_{S} c_1(L,h)\geq 0.
\] On the other hand,  the positivity of the extended metric allows to apply Chern--Weil type results, see e.g.\ \cite[Corollary 2.13]{BHdJ19} (and \cite[\S 5]{dJ4_ChernWeil} for the general framework), and we find 
\begin{equation}\label{positivity of height pairing}
\langle z, z^\vee\rangle_{X_K}=\deg(L)=\deg(L,{h})\geq 0.
\end{equation}  
\end{void}
\begin{rem}
	Let $z\in \CH^p_{\hom}(X/S)$. Let $c_1(\PP,||.||)$ be the translation-invariant curvature form of the biextension metric on the Poincaré bundle $\PP$ of $J^p_{X/S}$. Then
	\[
	\deg(L,{h})=\int_{S}(\nu_z,\nu_z^\vee)^*c_1(\PP,||.||).
	\] Hence, by \cite[Corollary 13.3]{Hain13}, $\deg(L,{h})=0$ if and only if the {Griffiths infinitesimal invariant}  $\delta(\nu_z)$ associated to the normal function $\nu_z\colon S\to J^p_{X/S}$ vanishes on $\Sbar$.
\end{rem}
\subsection{Conventions}

\begin{itemize}
	\item 	(Conventions for varieties) By an \textbf{algebraic scheme} over  $k$ (or an algebraic $k$-scheme), we mean a separated scheme of finite type over a field $k$. An \textbf{(algebraic) variety} over~$k$ is a geometrically integral algebraic scheme. A \textbf{nice variety} over $k$ is a geometrically connected smooth projective $k$-variety. An \textbf{(algebraic) curve} is a 1-dimensional variety. 
	
	\item (Complex varieties) By a $\C$-\textbf{variety}, we always mean an \emph{algebraic} variety over~$\C$ (as opposite to its associated complex-analytic structure). A \textbf{complex function field} is a field $k$ which is isomorphic to the function field $\C(S)$ of some smooth $\C$-variety~$S$.
	
	\item Let $S$ be a scheme resp.\ a complex manifold. By a \textbf{partial compactification} $\Sbar$ we mean that $\Sbar$ contains $S$ as a Zariski open resp.\ analytic open subset.  We will always assume  $\Sbar$ to be regular (but not necessarily proper). 
	
	\item We use the additive notation for Picard groups. For example, if $L$ and $M$ are line bundles, $L-M$ means $L\otimes M^\vee$.
	
	\item ($\Q$-line bundles) We will use the formalism of $\Q$-line bundles as explained in \cite[\S 2.1]{BHdJ18} or \cite[\S 2.2]{YZ21}. Following the above practice, we  think of and denote a $\Q$-line bundle somewhat informally in the form of $\frac{1}{m}L$, where $m\in \Z\ohne\{0\}$ and $L$ is  a usual line bundle.
\end{itemize}

\subsection*{Acknowledgements}
I would like to express my deep gratitude to Robin de Jong for suggesting the problem, for numerous inspiring conversations, and for his valuable remarks on an earlier version of this manuscript. I would like to thank Yinchong Song for his interest in Conjecture \ref{Conjecture: geometric height pairing} and for related discussions, see  \cite[\S 6]{Song2}. I wish to thank Klaus Künnemann for sharing the Diplomarbeit \cite{SeiboldM} and for his thoughtful comments on my presentation of the current work at the 2025 Intercity Seminar on Arakelov Geometry. I am also grateful to Alexandru Ioan Cuza University of Iași for its hospitality in hosting the seminar, which made many fruitful discussions possible. I would like to extend my thanks to David Holmes for clarifying several concepts appearing in the last section of this paper.  I acknowledge the financial support of the CSC--Leiden University Scholarship No. 202208080111.

\section{Algebraically trivial cycles revisit}\label{Algebraically trivial cycles revisit}
We recall the concept of correspondences and fix notations at the same time. The references are \cite[\S 16.1]{Fulton} and \cite[\S 1]{MNP}.

In the following let $X, Y, Z$ be  nice varieties  over a field. Thus we have the intersection product $\cdot$ on   Chow groups of these objects, and the proje\-ctions $\pr_{(-)}$ are proper and flat so that pushforward and pullback along these are well-defined.
\begin{defn}
	A \textbf{correspondence} $\alpha$ from $X$ to $Y$ is a cycle class $\alpha\in\CH^i(X\times Y)$. The \textbf{transpose} $\alpha^\top\in \CH^i(Y\times X)$  is given by the pushforward $\alpha^\top\coloneq \tau_*\alpha\in \CH^i(Y\times X)$, where 
	\[
	\tau\colon X\times Y\xrightarrow{\cong} Y\times X,\quad (x,y)\mapsto (y,x).
	\] Immediately by definition, one has 
	\[
	(\alpha^\top)^\top=\alpha.
	\]
\end{defn}
Every morphism $f\colon X\to Y$ (identified with its graph $\Gamma_f$) is a  correspondence from $X$ to $Y$. One can pullback and pushforward cycle classes along $f$. These operations work actually for more general correspondences:
\begin{void}[Action of correspondence]\label{Action of correspondence}
	Let $\alpha\in \CH^i(X\times Y)$ be a correspondence, it acts on cycles as follows.   For every $p\in \N$, we have the covariant action
	\begin{displaymath}\label{covariant action of correspondence}
		\alpha_*\colon\CH^p(X)\to \CH^{p+i-\dim X}(Y),\quad x\mapsto \pr_{Y,*}(\pr_X^*x\cdot \alpha).
	\end{displaymath} Dually, we have the contravariant action
	\begin{displaymath}
		\alpha^*\colon \CH^p(Y)\to \CH^{p+i-\dim Y}(X),\quad y\mapsto \pr_{X,*}(\pr^*_Yy\cdot \alpha).
	\end{displaymath}
	They are related by
	\begin{displaymath}\label{transpose and pullback}
		(\alpha^\top)_*=\alpha^* \quad\text{and dually}\quad \alpha_*=(\alpha^\top)^*.
	\end{displaymath}
\end{void}
\begin{void}[Composition of correspondence]\label{composition of correspondence}
	Suppose $\alpha\in \CH^i(X\times Y), \beta\in \CH^j(Y\times Z)$ are correspondences. They can be composed by the rule
	\[
	\beta\circ \alpha\coloneq \pr_{XZ,*}(\pr_{XY}^*\alpha\cdot \pr_{YZ}^*\beta)\in \CH^{i+j-\dim Y}(X\times Z),
	\] where $\pr_{(-)}$ denotes the natural projection, e.g.\ $\pr_{XY}\colon X\times Y\times Z\to X\times Y$.
\end{void}

Next, we recall the notion of algebraically trivial cycles and fix our definition at the same time. Indeed, there have been different conventions for this notion across the literature, which are not always equivalent to each other (when the ground field is not algebraically closed, say). As elaborated in \cite{ACV19Alg}, the situation  is more favourable  if the ground field is perfect. 

In the rest  of this section, $k$ denotes a perfect field and $k\subset K$ a fixed algebraic closure.
\begin{void}[Refined Gysin pullback]\label{Refined Gysin pullback}
	Let $T$ be a \emph{smooth} $k$-variety of dimension~$e$, then each $k$-point $t\in T(k)$ induced a closed regular embedding $t\colon \Spec k\hookrightarrow T$ of codimension $e$. Let $X$ be an algebraic $k$-scheme and consider the  pullback diagram
	\[
	\begin{tikzcd}
		X_t\cong X=(T\times_k X)\times_T \Spec k\ar{d}\ar{r} & T\times_k X\ar{d}\\
		\Spec K\ar{r}{t} & T
	\end{tikzcd}
	\] where $X_t$ denotes the fiber of $T\times_K X$ over $t$. By \cite[\S 6.2]{Fulton}, we have  a well-defined homomorphism
	\[
	t^!\colon Z_i(T\times_K X)\to \CH_{i-e}(X)
	\] which descends to
	\[
	t^!\colon \CH_i(T\times_K X)\to \CH_{i-e}(X),
	\] called the \textbf{refined Gysin pullback}.  For $\alpha\in \CH_i(T\times_K X)$, we also write $t^!\alpha=\alpha_t$. 
\end{void}
\begin{defn}\label{~alg}
	Let $X$ be an algebraic $k$-scheme. A cycle class $z\in \CH_i(X)$ is called \textbf{algebraically trivial} over $k$, if there is a nice $k$-\emph{curve} $T$, $k$-points $t_0,t_1\in T(k)$ and $\alpha\in \CH_{i+1}(T\times_K X)$ such that
	\[
	z=\alpha_{t_1}-\alpha_{t_0}.
	\]
	A cycle $Z\in \Zykle_i(X)$ is called \textbf{algebraically trivial} over $k$, if its associated cycle class $z=[Z]$ is algebraically trivial over $k$. The subgroups of algebraically trivial cycles resp.\ cycle classes will be denoted by $\Zykle_{i,\alg}(X)\subset \Zykle_i(X)$ resp.\ $\CH_{i,\alg}(X)\subset \CH_i(X)$. We introduce $\Zykle^i_\alg(X)$ and $\CH^i_\alg(X)$ in the same way.
	
See \cite[Proposition 3.14]{ACV19Alg} for more equivalent definitions of the algebraic triviality using different parameter spaces $T$.
\end{defn}

\begin{prop}\label{Lemma 3.18, ACV}
	Let $X$ be an algebraic $k$-scheme and $Z\in \Zykle_{i}(X)$ be a cycle on $X$. The following properties for $Z$ are equivalent:
	\begin{enumerate}
		\item  $Z$ is algebraically trivial over $k$.
		\item There is a nice $k$-curve $T$, an effective  cycle $\mathscr{Z}\in \Zykle_{i+1}(T\times_K X)$ which is flat over $T$ and $t_0,t_1\in T(k)$ such that
		\[
		Z=\mathscr{Z}_{t_1}-\mathscr{Z}_{t_0}.
		\] (Here $\mathscr{Z}_{t_i}$ denote the usual scheme-theoretic fibers.)
	\end{enumerate}
\end{prop}
\begin{proof}
	See \cite[Lemma 3.18]{ACV19Alg}.
\end{proof}
\begin{rem}\label{alg. trivial is adequate}
	Algebraic triviality of algebraic cycle resp.\ algebraic cycle classes on an algebraic $k$-scheme $X$ defines an adequate equivalence relation (see e.g.\ \cite[\S 1.2]{MNP}) on the graded groups $Z(X)=\bigoplus_i\Zykle^i(X)$ resp.\  $\CH(X)=\bigoplus_i \CH^i(X)$. In particular, it is preserved under correspondences. 
\end{rem}
\begin{conv}
	When we say that a cycle (class) is algebraically trivial on an  algebraic scheme $X$ without explicit reference of the field, it is assumed to be algebraically trivial over the field of definition of $X$.
\end{conv}
\begin{defn}\label{Relative algebraic triviality, general}
	Let $X\to S$ be a proper faithfully flat morphism over a \emph{smooth} $k$-variety~$S$. A relative algebraic cycle $Z\in \Zykle^\bullet(X/S)$ is called \textbf{(fiberwise or relative) algebraically trivial}, if for every geometric point $t\in S(K)$, the  fiber~$Z_t$ is algebraically trivial  on~$X_t$. Denote the resulting subgroup by $\Zykle^\bullet_\alg(X/S)\subset \Zykle^\bullet(X/S)$. A relative cycle class $z\in \CH^\bullet(X/S)$ is  \textbf{(fiberwise or relative) algebraically trivial}, if its Gysin pullback $z_t=t^!z\in \CH^\bullet_\alg(X_t)$ for every geometric point $t\in S(K)$. Denote the resulting subgroup by $\CH^\bullet_\alg(X/S)\subset \CH^\bullet(X/S)$.
\end{defn}
\begin{rem}[Gysin fiber vs.~scheme fiber]\label{repn of Gysin pullback}
	Suppose $z\in \CH^\bullet(X/S)$ is represented by a relative cycle $Z\in \Zykle^\bullet(X/S)$. Then by the smoothness of $S$, for every geometric point $t\in S(K)$, the  class $[Z_t]$ of the fiber $Z_t$ represents the Gysin pullback $z_t$, see \cite[p.176]{Fulton} or the paragraph between Definitions 3.16 and 3.17 of \cite{ACV19Alg}. Consequently, a cycle class $z\in \CH^\bullet_\alg(X/S)$ if and only if it can be represented by a relative cycle $Z\in\Zykle^\bullet_\alg(X/S)$.
\end{rem}

\begin{rem}\label{Weil-Bloch}
	Let $X$ be a nice $k$-variety.  Proposition \ref{Lemma 3.18, ACV} and Remark \ref{repn of Gysin pullback} imply that for $z\in \CH_{i,\alg}(X)$, we can find a nice $k$-curve~$T$, a correspondence $\alpha\in \CH_{i+1}(T\times_k X)$ and $t_0, t_1\in T(k)$ such that 
	\[
	z=t_1^!\alpha-t_0^!\alpha=\alpha_*([t_1]-[t_0])=(\alpha^\top)^*([t_1]-[t_0]).
	\] 
\end{rem}
The  next lemma roughly says that one can spread out relative algebraic triviality. 
\begin{lem}\label{spread out alg. triv.}
	Let  $X\to S$ be a proper faithfully flat morphism over a smooth $k$-variety~$S$. The following conditions are equivalent for a relative cycle $Z\in \Zykle^p(X/S)$:
	\begin{enumerate}
		\item For every geometric point $t\in S(K)$, the geometric fiber $Z_t$ is algebraically trivial.
		\item For some geometric point $t\in S({K})$, the geometric fiber $Z_t$ is algebraically trivial.
		\item The geometric generic fiber $Z_{K}$ is algebraically trivial.
	\end{enumerate}
\end{lem}
\begin{proof}
	This is a special case of \cite[Lemma 1.5]{ACV23}:  We substitute $\Lambda=k$, $S=\Spec k$ and $T=S$, where the left hand side are the notations  in the cited lemma and right hand side  the notations of the present lemma, and we consider the family $X\times_K S\to S$ (in the notation of the present lemma).
\end{proof}
\section{Regular homomorphisms and higher Picard varieties}\label{Regular homomorphisms and higher Picard varieties}
 The concept of regular homomorphisms were initiated by Weil and Griffiths as a tool to analyse the structure of  algebraic cycles. A  systematic investigation over $\C$ was first found in Lieberman's work \cite{Lieberman} (see also the survey \cite{LiebermanSurvey}). Later,  H. Saito gave a purely algebraic approach in \cite{HSaito}; we note that M. Seibold has worked out many details of the latter in his Diplomarbeit \cite{SeiboldM}. The concept has been recently revisited from a more functorial point of view in the framework  \cite{ACV23} of Achter--Casalaina-Martin--Vial. We will also review the notion of higher Picard varieties, which are representing objects of a subclass of regular homomorphisms, called Picard homomorphisms (Definition \ref{Def. Picard variety}). Up to incidence equivalence, the latter varieties parametrize algebraically trivial cycles (Lemma \ref{higher Picard and algebraically trivial cycles}). An important application is that one can define the arithmetic Beilinson--Bloch height pairing  for algebraically trivial cycles unconditionally, a theorem due to Künnemann (\cite{Kue96}), see also \cite[Remark 4.0.8]{BeilinsonHP}. 
\subsection{Relative setting}
\begin{notn}\label{Notation for ACV's functor}
	\hfil
	\begin{enumerate}
		\item[(1)] Unless stated differently, in this section $k$ denotes a perfect field, $S$ a smooth $k$-variety.
		\item[(2)] Let $\mathsf{Ab}$ be the category of abelian groups. Let $\mathsf{Sm}/S\coloneq \mathsf{Sm}_k/S$ be the category whose objects are separated $S$-sche\-mes that are smooth and of finite type over $k$ with structure maps to $S$ being dominant. Morphisms in $\mathsf{Sm}/S$ are $S$-morphisms.
		\item[(3)] Given  $S$-schemes $X\to S$ and $T\to S$, we write $X_T$ for $T\times_S X$.
	\end{enumerate}
\end{notn}
Recall the notation $\CH^p_\alg(X_T/T)$   ($T\in \mathrm{Sm}/S$) from Definition \ref{Relative algebraic triviality, general}.
\begin{defn}\label{Functor of algebraically trivial cycles}
	Let  $X\to S$ be a proper faithfully flat morphism of relative dimension~$d$.   Let $p\in\{1,\ldots,d\}$. The \textbf{functor of codimension}-$p$ \textbf{algebraically trivial cycles} on the family $X/S$ is the \emph{contra\-variant} functor
	\[
	\A^p_{X/S}\colon \mathsf{Sm}/S\to \mathsf{Ab}, \quad T\mapsto \CH^p_\alg(X_T/T).
	\] A morphism $f\colon T'\to T$ in $\mathsf{Sm}/S$ is to be sent to the refined Gysin pullback
	\[
	f^!\colon \CH^p_{\alg}(X_T/T)\to \CH^p_\alg(X_{T'}/T').
	\]  
\end{defn}
\begin{defn}\label{Regular homomorphism, functorial}
	Let $A/S$ be an abelian scheme over $S$, viewed as a contravariant functor  $\mathsf{Sm}/S\to \mathsf{Ab}$. Let $X/S$ be a proper faithfully flat family as before. A \textbf{regular homomorphism} (in codimension $p$) on $X/S$ is a morphism of functors $$\Phi^p\colon \A^p_{X/S}\to A/S.$$ An \textbf{algebraic representative} of $\A^p_{X/S}$ in codimension $p$ is an abelian scheme $\text{Ab}^p_{X/S}$ together with a morphism of functors
	\[
	\Phi^p_{X/S}\colon  \A^p_{X/S}\to \text{Ab}^p_{X/S}
	\] which is initial among all morphisms of functors $\Phi^p\colon \A^p_{X/S}\to A/S$. That is, for all such $\Phi^p$ there exists a unique $S$-morphism of abelian schemes $f\colon \mathrm{Ab}^p_{X/S}\to A$ such that the diagram commutes:
	\[
	\begin{tikzcd}
		\A^p_{X/S}\ar{r}{\Phi^p_{X/S}} \ar[swap]{d}{\forall\,\Phi^p} & \text{Ab}^p_{X/S}\ar[dashed]{ld}{\exists ! f}\\
		A/S
	\end{tikzcd}
	\]
\end{defn}
Since algebraic triviality is preserved under correspondences, using \ref{Action of correspondence}, one immediately sees the following
\begin{lem}\label{regular homomorphism and correspondence}
	Let $X, T$ be   nice varietyies over $k$ and $\alpha\in \CH^i(X\times_k T)$ be a correspondence. There is a morphism of functors
	\[
	\A^{p+\dim T-i}_{T/k}\to \A^p_{X/k}
	\] induced by $\alpha^*$.  Similarly, there is a morphism of functors
	\[
	\A^{p}_{X/k}\to \A^{p+i-\dim X}_{T/k}
	\] induced by $\alpha_*$.  (For all $p$ such that the functors are defined.)
\end{lem}

\subsection{Classical setting}
\hfil

For later  purposes, let us specify the notion to the classical situation where $S=\Spec K$ with an algebraically closed field $K$. The functor of algebraically trivial cycles will be simply denoted by $\A^\bullet_X$. 
\begin{defn}\label{Regular homomorphism, classical}
	Suppose we are given  an abelian $K$-variety $A$,  a nice $K$-variety $X$ of dimension $d$.  A group homomorphism 
	\[
	h\colon \CH^p_\alg(X)\to A(K),
	\] is said to be a  \textbf{regular homomorphism} (in codimension $p\in\{1,\ldots,d\}$) if for all smooth  $K$-varieties $T$, all $K$-points $t_0\in T(K)$ and all correspondences $\alpha\in \CH^p(T\times_K X)$, the map
	\begin{equation}\label{test regular hom}
		\varphi\coloneq \varphi_{T,\alpha,t_0}\colon	T(K)\to \CH^p_\alg(X)\xrightarrow{h} A(K),\quad t\mapsto h(\alpha_*([t]-[t_0]))=h(\alpha_t-\alpha_{t_0})
	\end{equation}
	is induced by a morphism $f\colon T\to A$ of $K$-schemes, i.e. $f(K)=\varphi$. Note that there can exist at most one such morphism.
\end{defn}
\begin{lem}\label{Regular homomorphisms, Basepoint-free definition}
	Suppose we are given  an abelian $K$-variety $A$,  a nice $K$-variety $X$ and a group homomorphism 
	\[
	h\colon \CH^p_\alg(X)\to A(K).
	\] The following claims are equivalent.
	\begin{enumerate}
		\item $h$ is a regular homomorphism in the sense of Definition \ref{Regular homomorphism, classical}.
		\item (Basepoint-free version) For all smooth  $K$-varieties $T$  and all correspondences $\alpha\in\A^p_X(T)=\CH^p_\alg(T\times_K X)$, the map 
		\[
		\varphi=\varphi_{T,\alpha}\colon T(K)\to \CH^p_\alg(X)\xrightarrow{h} A(K),\quad t\mapsto h(\alpha_t)
		\] is induced by a morphism $T\to A$ of $K$-schemes.
	\end{enumerate}
\end{lem}
\begin{proof}
	See \cite[Lemma 1.12]{ACV23}.
\end{proof}
\begin{rem}\label{translation of functorial regular hom}
	Lemma \ref{Regular homomorphisms, Basepoint-free definition} implies that there is  a morphism of functors $\Phi^p\colon \A^p_X\to A$ such that $h=\Phi^p_X(K)$. Indeed, we can define the morphism by the rule
	\[
	\Phi^p(T)(\alpha)\coloneq \varphi_{T,\alpha}\in A(T)=\mathrm{Hom}_K(T,A),
	\] using the notations of Lemma \ref{Regular homomorphisms, Basepoint-free definition}.
\end{rem}

\begin{rem}
	It is natural to ask whether the functor $\A^p_X$ admits algebraic representatives. A main result of Achter--Casalaina-Martin--Vial, \cite[Theorem 2]{ACV23}, implies that for a smooth projective family $X/S$ over a smooth $K$-variety $S$, where $K$ is a perfect field, the functor $\A^p_{X/S}$ admits an algebraic representative for $p=1,2$ and $d$, where $d$ is the relative dimension of  $X/S$. 
\end{rem}
\subsection{Higher Picard varieties}
\hfil

Let us go back to the classical case where $S=\Spec K$, $K$ is an \emph{algebraically closed} field and~$X$ is a nice  $d$-dimensional $K$-variety. In this section we will see that the functor $\A^p_X$ admits an algebraic representative for  $p\in \{1,\ldots,d\}$ in a weaker sense, namely when the test regular homomorphisms are replaced by the so-called Picard homomorphisms. 
\begin{defn}\label{def. Picard homomophism}
	Let $A$ be an abelian variety over $K$. A group homomorphism 
	\[
	h\colon \CH^p_\alg(X)\to A(K).
	\]   Then
	$h$ is said to be a \textbf{Picard homomorphism} if there exists a nice  $K$-variety $Y$, a correspondence $\beta\in \CH^{d+1-p}(X\times_K Y)$ and a closed immersion $$A\hookrightarrow \Pic^0_{Y/K,\mathrm{red}}$$ into the Picard variety  of $Y$ such that  the following diagram commutes:
	\begin{equation}\label{commutativity of Picard hom}
		\begin{tikzcd}
			\CH^p_{\alg}(X)\ar{d}{h} \ar{r}{\beta_*} & \CH^1_\alg(Y)\ar{d}{\cong}\\
			A(K)\ar[hookrightarrow]{r} & \Pic^0_{Y,\mathrm{red}}(K)
		\end{tikzcd}
	\end{equation}
\end{defn}
	\begin{defn}\label{Def. Picard variety}
A $p$-th \textbf{Picard variety} is an abelian variety $P$ over $K$ together with a surjective Picard homomorphism $\theta^p\colon \CH^p_\alg(X)\twoheadrightarrow P(K)$ that is \emph{initial} for  all Picard homomorphisms, i.e.\ for all Picard homomorphisms $h\colon \CH^p_\alg(X)\to A(K)$, there is a unique $K$-morphism of abelian varieties $f\colon P\to A$ such that the diagram commutes:
		\[
		\begin{tikzcd}
			\CH^p_\alg(X)\ar[swap]{d}{\forall\, h} \ar{r}{\theta^p} & P(K) \ar[dashed]{ld}{\exists ! f(K)}\\
			A(K) & 
		\end{tikzcd}.
		\] A \textbf{higher Picard variety} is a $p$-th Picard variety for some $p\in \{1,\ldots, d\}$, denoted by $\Pic^p_{X/K}$ or simply by $\Pic^p_X$. The universal map $\theta^p$ will be referred as the \textbf{algebraic Abel--Jacobi map} (in codimension $p$).
	\end{defn}
	\begin{lem}\label{Every Picard is regular}
		Every Picard homomorphism is a regular homomorphism.
	\end{lem}
	\begin{proof}
		Let $h\colon \CH^p_\alg(X)\to A(K)$ be a Picard homomorphism. Let $\varphi=\varphi_{T,\alpha,t_0}$ be a test map, notated as in Definition \ref{Regular homomorphism, classical}. We want to show that the map
		\[
		\varphi\colon T(K)\to A(K)
		\] is induced by a morphism $f\colon T\to A$. It is enough to show that the composite
		\[
		T(K)\to A(K)\hookrightarrow \Pic^0_{Y,\mathrm{red}}(K)
		\] is induced by a morphism $f'\colon T\to \Pic^0_{Y,\mathrm{red}}$ of schemes, since $T$ is reduced and the set-theoretic image $f'(T)$ is contained in $A$, \emph{cf}.\ \cite[Remark 10.32]{GW}. 
		
	Let $Y$ and $\beta$ be as in Definition \ref{Def. Picard variety}.	Using the commutativity of \eqref{commutativity of Picard hom}, it is equivalent to showing that
		\[
		\varphi'\colon T(K)\to \CH^p_\alg(X)\xrightarrow{\beta_*} \CH^1_\alg(Y)\cong \Pic^0_{Y,\mathrm{red}}(K)
		\] is induced by a morphism of schemes. we know $\CH^1_\alg(Y)\cong \Pic^0_{Y,\mathrm{red}}(K)$ is a regular homomorphism, i.e. there is a morphism of functors
		\[
		\A^1_Y\to \Pic^0_{Y,\mathrm{red}}
		\] by Remark \ref{translation of functorial regular hom}.
		 Applying Lemma \ref{regular homomorphism and correspondence} to the composite $\beta\circ \alpha\in \CH^1(T\times Y)$, we obtain a morphism of functors
		 \[
		 \A^{\dim T}_{T}\xrightarrow{(\beta\circ \alpha)_*} \A^1_Y\to  \Pic^0_{Y,\mathrm{red}}.
		 \] This implies that $\varphi'$ is induced by a morphism.
	\end{proof}
\begin{prop}
	Higher Picard varieties exist in every codimension $p\in \{1,\ldots d\}$.
\end{prop}
\begin{proof}[Comments]
The basic idea is to use Saito's criterion, see \cite[Theorem 3.6]{HSaito} or \cite[Proposition 5.3]{ACV23}. 

First, for every  regular homomorphism $h\colon \CH^p_\alg(X)\to A(K)$ there is a unique abelian subvariety $C\subset A$ over $K$ with the property that  $h$ factors through $B(K)$:
\[
h=(\CH^p_\alg(X)\twoheadrightarrow B(K) \stackrel{\iota(K)}{\hookrightarrow}A(K)).
\] If $K=\C$, this is outlined in \cite[Proposition 12.23]{LewisBook}. The general case is elaborated in \cite[\S V, Lemma 1.3]{SeiboldM}. We define $\dim \mathrm{im}(h)\coloneq \dim B$.

Choose a Weil cohomology theory $H(-)$ over a field $F$ of characteristic zero. One can show that 
\begin{equation}
	\dim \mathrm{im}(h)\leq \frac{1}{2}\dim_F H^{2p-1}(X)
\end{equation} for every \emph{Picard homomorphism} $h$. Saito's criterion then implies the existence of the representing object, namely the higher Picard variety.
\end{proof}

\begin{ex}[Higher Picard varieties interpolate Albanese and Picard varieties.]\label{Picard and higher Picard} 
By construction, we have $\Pic^d_{X/K}=\mathrm{Alb}_{X/K}$ resp.\ $\Pic^1_{X/K}=\Pic^0_{X/K,\mathrm{red}}$ the Albanese resp.\ Picard variety associated to $X$, where  $(-)^0$ denotes the identity connected component of a group scheme. And the algebraic Abel--Jacobi maps in codimension $d$ resp.\ 1 are the usual universal maps. We accept a slight notational conflict: $\Pic^1_{X/K}=\Pic^0_{X/K,\mathrm{red}}$.  
\end{ex}
	
\begin{void}[Functoriality]\label{functoriality of higher Picard varieties}
		Let $X, Y$ be nice $K$-varieties of dimension $d, e$ respectively. Let $\alpha\in \CH^{e-p+q}(Y\times X)$ be a correspondence. By the universal property of higher Picard varieties, there is a unique morphism $\Pic(\alpha)\colon \Pic^p_{Y}\to \Pic^q_{X}$ (for $p\in \{1,\ldots, e\},\, q\in \{1,\ldots, d\}$) making the following diagram commute
		\begin{equation}\label{functoriality of Pic()}
			\begin{tikzcd}
				\CH^p_\alg(Y)\ar{rr}{\alpha_*}\ar{d}{\theta^p} & & \CH^q_\alg(X)\ar{d}{\theta^q}\\
				\Pic^p_{Y}(K)\ar{rr}{\Pic(\alpha)}& & \Pic^q_{X}(K)
			\end{tikzcd}
		\end{equation}
\end{void}
	\begin{void}
		Putting $p=e$, the above functoriality yields a  group homomorphism
		\[
		\CH^q(Y\times X)\to \mathrm{Hom}(\mathrm{Alb}_{Y},\Pic^q_{X}).
		\] Specializing  further to $Y=\Pic^p_{X}$ and $q=p$, and noting $\mathrm{Alb}_{\Pic^p_X}=\Pic^p_X$, we obtain a homomorphism
		\[
		\beta\colon \CH^p(\Pic^p_{X}\times X)\to \mathrm{End}(\Pic^p_{X}).
		\]
	\end{void}
	\begin{lem}\label{k^p_X}
		There is a positive integer $k^p_X$ such that
		\[
		\mathrm{im}(\beta)\cap \langle \id_{\Pic^p_{X/K}}\rangle=\langle k^p_X\rangle,
		\] where $\langle,\rangle$ denotes the $\Z$-module generated by the respective elements as a subgroup of $\mathrm{End}(\Pic^p_{X})$. 
	\end{lem}
	\begin{proof}
		See \cite[(4.1)]{HSaito}.
	\end{proof}
	\begin{defn}\label{Poincare p-cycle}
		Any element $\alpha_X^p\in \CH^p(\Pic^p_{X}\times X)$ that is mapped to  $k_X^p\cdot 1\in \mathrm{End}(\Pic^p_X)$ under $\beta$  is called a \textbf{Poincaré} $p$-\textbf{cycle}. 
	\end{defn}
	\begin{void}\label{Pic of Poincare p-cycle}
		Let $\alpha_X^p$ be a Poincaré $p$-cycle.  Let $d\coloneq \dim X,\, e\coloneq \dim \Pic^p_{X}$. By definition, we have a commutative diagram ($[n]$ denotes the multiplication-by-$n$ map) 
		\begin{equation}
			\begin{tikzcd}
				\CH^e_\alg(\Pic^p_{X})\ar{rr}{(\alpha_X^p)_*}\ar{d}{\theta^e} & & \CH^p_\alg(X)\ar{d}{\theta^p}\\
				\Pic^p_{X}(K)\ar{rr}{[k_X^p]}& & \Pic^p_{X}(K)
			\end{tikzcd}
		\end{equation}
		Consider $q\coloneq d+1-p$, $(\alpha^p_X)^\top\in \CH^p(X\times \Pic^p_{X})$ and $\lambda^q_X\coloneq \Pic((\alpha^p_X)^\top)$. There is  a commutative diagram
		\begin{equation}\label{Pic(alpha^*)}
			\begin{tikzcd}
				\CH^q_\alg(X)\ar{rr}{(\alpha_X^p)^*}\ar{d}{\theta^q} & & \CH^1_\alg(\Pic^p_{X})\arrow[d, "\cong"', "\theta^1"]\\
				\Pic^q_{X}(K)\ar{rr}{\lambda_X^q}& & \Pic^1_{\Pic^p_{X}}(K)=\Pic^{p,\vee}_{X}(K)
			\end{tikzcd}
		\end{equation}
	\end{void}
	\begin{lem}\label{properties of Poincare cycles}
		Let $X, Y$ be nice $K$-varieties of dimension $d, e$. We have 
		
		\begin{enumerate}
			\item Every choice of the Poincaré $p$-cycle $\alpha^p_X$ induces the same  morphism $$\lambda_X^{d+1-p}\colon \Pic^{d+1-p}_{X}\to (\Pic^p_{X})^\vee.$$
			\item  $\lambda_X^p$ is an  isogeny (it is an isomorphism if $p=d$).
			\item Given a correspondence $\alpha\in \CH^{p+q}(Y\times X)$, one has
			\[
			[k_X^p]\circ \lambda_Y^{q+1}\circ\Pic(\alpha^\top)=[k_Y^{e-q}]\circ \Pic(\alpha)^\vee\circ \lambda_X^{d+1-p},
			\]  as morphisms  $\Pic^{d+1-p}_{X}\to \Pic^{e-q}_{Y}$
		\end{enumerate} (for all $p,q\in \N$ such that the notations make sense).
	\end{lem}
	\begin{proof}
		See  \cite[\S 4]{HSaito}.
	\end{proof}
	\begin{ex}\label{Polarization on curves}
		Let $C$ be a nice $K$-curve and $\lambda_C\colon \Pic^1_{C}\xrightarrow{\cong}(\Pic^1_{C})^\vee$ be the  principal polarization induced by the theta divisor. Then $k_C^1=1$ and $\lambda^1_C=[-1]\circ \lambda_C$.
	\end{ex}
	\begin{rem}[Higher Picard varieties over a perfect field]\label{Higher Picard var over perfect field}
		Let $X$ be a nice variety over a perfect but not necessarily algebraically closed field $k$ (e.g.\ a complex function field); let $k\subset K$ be an algebraic closure. Using  the language of Galois descents, one can construct an abelian variety $\Pic^p_{X/k}$ associated to  $X$ such that
		\[
		\Pic^p_{X/k}\otimes_k K\cong \Pic^p_{X_K/K}
		\] is the $p$-th Picard variety for every $p\in\{1,\ldots, \dim X\}$ introduced before; and the construction is $\mathrm{Gal}(K/k)$-equivariant in some sense. The details has been worked out in the last chapter (\S V) of Seibold's Diplomarbeit \cite{SeiboldM}.
	\end{rem}
\section{Relative algebraic intermediate Jacobians}\label{Algebraic intermediate Jacobians}
The goal of this section is to explain the fact that Abel--Jacobi trivial cycles are incidentally trivial in the complex function field setting. This will yield a natural surjective morphism
\[
h^p\colon J^p_{a,X}\twoheadrightarrow \Pic^p_{X}
\] from the geometric generic fiber of the relative algebraic intermediate Jacobian to the higher Picard variety. Here $X$ denotes  a nice variety over an algebraically closed field $K$ of a complex function field, and we fix an embedding $K\subset \C$.  The similar statement is true when the base field is $\C$, and the claim follows  from a descent argument from this classical fact by descent arguments. 

Griffiths's Conjecture \ref{Griffiths conjecutre for incidentally trivial cycles} implies that $h^p$ is an isogeny, which turns out to be necessary to complete the final comparison.

We begin by recalling the domain of $h^p$ based on the theory of relative algebraic intermediate Jacobians developed by Achter–Casalaina-Martin–Vial.
\begin{void}[Abel--Jacobi map]
	Let $X$ be a projective complex manifold of dimension $d$. Fix $p\in\{1,\ldots,d\}$ and consider the (singular) cohomology group $H^{2p-1}(X,\Z)$. The Hodge decompo\-sition provides a Hodge (decreasing) filtration $F^\bullet$ on $$H^{2p-1}(X,\C)\cong H^{2p-1}(X,\Z)\otimes \C.$$ 	The $p$-th \textbf{Griffiths intermediate Jacobian} is given by
	\[
	J^p_X\coloneq H^{2p-1}(X,\C)/(F^p\oplus H^{2p-1}(X,\Z)_\mathrm{tf}),
	\] where $F^p\coloneq F^p H^{2p-1}(X,\C)$  and  $A_\mathrm{tf}$ denotes the torsion-free part of an abelian group $A$. The quotient is in fact a compact (albeit not necessarily algebraic) complex torus  that admits the \textbf{Abel--Jacobi map}: 
	\[
	\AJ^p\colon \Zykle^p_{\hom}(X)\to J^p_X,
	\] where $Z^p_{\hom}(X)$  denotes the subgroup of homologically trivial (algebraic) cycles. It can be shown that $\AJ^p$ factors through $\CH^p_{\hom}(X)$, and the induced map will still be denoted by $\AJ^p$.  See Murre's lecture \cite[\S 9.3]{MurreLecture} for more details.
\end{void}
\begin{void}[Normal function associated to an algebraic cycle]
	Placing the Abel–Jacobi map in a family yields a main source of (admissible) normal functions. More precisely, consider a projective submersive (holomorphic) map $f\colon X\to S$ over a complex manifold $S$; let $d$ be the fiber dimension. Fix $p\in \{1,\ldots,d\}$. The family $X/S$ gives rise to  a polarized $\Z$-variation of Hodge structure with underlying local system $R^{2p-1}f_*\Z_X$, its fiber over $s\in S$ is exactly $H^{2p-1}(X_s,\Z)$ and one can form the relative torus $J^p_{X/S}\to S$ such that its fiber $J^p_{X/S,s}=J^p_{X_s}$ is the $p$-th Griffiths intermediate Jacobian of the fiber $X_s$. By a classical result of Griffiths, the structure map $J^p_{X/S}\to S$ is  holomorphic.  For $z\in \CH^p_{\hom}(X/S)$, we can consider its normal function
	\[
	\nu_z\colon S\to J^p_{X/S},\quad s\mapsto \AJ^p_s(z_s),
	\] where $\AJ^p_s\colon \CH^p_{\hom}(X_s)\to J^p_{X_s}$ is the Abel--Jacobi map considered above. Again by Griffiths, $\nu_z$ is a holomorphic section of the relative torus $J^p_{X/S}\to S$; the map $\nu_z$ is  an \textbf{admissible normal function}, i.e.\ it corresponds to a certain admissible variation of mixed Hodge structures, by \cite[Remark 1.7]{Saito}.
\end{void}
Recall that every algebraically trivial algebraic cycle is homologically trivial. The study of the Abel–Jacobi map restricted to such cycles is natural and interesting, and has been carried out in a series of works by Achter–Casalaina-Martin–Vial. Below we present (and slightly reformulate) some of their main results. 

Let $(-)^\an$ denote the complex-analytification functor.
\begin{thm}[Achter--Casalaina-Martin--Vial]\label{ACV, Main Theorems}
	Let $f\colon X\to S$ be a smooth projective morphism of relative dimension $d$ over a smooth $\C$-variety $S$. Let $p\in\{1,\ldots,d\}$. There is a unique  abelian scheme $J^p_{a,X/S}\to S$ with the following properties.
	\begin{enumerate}
		\item	The holomorphic map $(J^p_{a,X/S})^\an\to S^\an$ is a relative subtorus of  $J^p_{X^\an/S^\an}\to S^\an$. That is, over each $s\in S^\an$, the fiber $(J^p_{a,X/S})^\an_s\subset (J^p_{X^\an/S^\an})_s=J^p_{X^\an_s}$ is a complex subtorus. Moreover, the construction of $J^p_{a,X/S}$ is compatible with base change $T\to S$ for $T\in \mathsf{Sm}/S$ (Notation \ref{Notation for ACV's functor}(2)), i.e.
		\[
		J^p_{a,X/S}\times_S T=J^p_{a,X_T/T}.
		\]
		\item The Abel--Jacobi map on algebraically trivial cycles is a surjective regular homo\-morphism, i.e.  there is a unique morphism of functors \[
		\Phi_{\AJ}^p\colon \A^p_{X/S}\to J^p_{a,X/S}
		\] which is surjective on geometric points (\cite[Definition 4.1]{ACV23}) such that for  $T\in \mathsf{Sm/S}$ and $z\in \CH^p_\alg(X_T/T)$, the holomorphic normal function
		\[
		\nu_z\colon T^\an\to J^p_{X_T^\an/T^\an}, \quad t\mapsto \AJ_t(z_t)
		\] factors through $(J^p_{a,X_T/T})^\an\hookrightarrow J^p_{X_T^\an/T^\an}$ and  is algebraic.  More precisely, there is an algebraic section $\nu_{z,a}\coloneq\Phi^p_\AJ(T)(z)\in J^p_{a,X_T/T}(T)$ of the abelian scheme $J^p_{a,X_T/T}$ such that
		\[
		\nu_z=\left(T^\an\xrightarrow{\nu_{z,a}^\an} (J^p_{a,X_T/T})^\an\subset J^p_{X_T^\an/T^\an}\right).
		\]	
		
		\item  The abelian scheme $J^p_{a,X/S}\to S$ is the algebraic representative of $\A^p_{X/S}$ for $p=1,2,$ and $d$.
	\end{enumerate}
\end{thm}
\begin{proof}
	See \cite[Theorems 9.1, 9.3 and 9.5]{ACV23}.
\end{proof}
\begin{defn}
	We call $\nu_{z,a}\coloneq\Phi^p_\AJ(T)(z)\in J^p_{a,X_T/T}(T)$ the \textbf{algebraic normal function} associated to $z$ (relative to $T\in \mathsf{Sm}/S$). 
\end{defn}
\begin{rem}
	\hfil
	\begin{enumerate}
		\item[(1)] By a classical theorem of Griffiths (see \cite[Corollary 12.25]{LewisBook}), the Abel--Jacobi map  on relative homo\-logically trivial cycles is a regular homomorphism in the holomorphic  category. Theorem \ref{ACV, Main Theorems}(2) provides an algebraic analogue.
		\item[(2)] In this item let $X$ be a single nice $\C$-variety. By our definition, $\AJ^p$ maps $\CH^p_\alg(X)$ surjectively onto $J^p_{a,X}(\C)$, and the image has the structure of a complex abelian variety. In the literature, ``algebraic intermediate Jacobian'' sometimes refers to the \emph{maximal} abelian variety contained in $J^p_{X}$; we denote this object temporarily by $J^p_{\mathrm{alg},X}$. By definition, $J^p_{a,X}\subset J^p_{\mathrm{alg},X}$, and the equality  follows from the Hodge Conjecture (see \cite[\S 12.2]{VoisinHodgeI}).
		\item[(3)] Let again $X/S$ be a smooth projective family. Unlike Griffiths intermediate Jacobians, the identification
		\[
		(J^p_{a,X/S})_s=J^p_{a,X_s}
		\] (the left hand side denotes the fiber of the abelian scheme $J^p_{a,X/S}$)  holds only for  ``very general'' points $s\in S(\C)$, see \cite[Remark 5.2]{ACV19} for the details. 
	\end{enumerate} 
\end{rem}
In the rest of this section we use the following
\begin{notn}\label{notation for fields for algebraic intermediate Jacobian}
	Let $k$ be the function field of a smooth $\C$-variety $S$ and fix an algebraic closure $K$ of $k$ inside $\C$. If  $X/S$ is a smooth projective family, we use $X_k$ resp.\ $X_K$ to denote the generic resp.\ geometric generic fiber.
\end{notn}
\begin{void}[Distinguished model]\label{distinguished model}
	Clearly, once we have the abelian scheme  $J^p_{a,X/S}/S$, we can consider its generic fiber $J^p_{a,X_k}$. In reality, the construction happens the other way around. To fix ideas, let $X$ be a nice variety over $k$. Let $J^p_{a,X_\C}$ be the classical algebraic intermediate Jacobian. By a main result of Achter--Casalaina-Martin--Vial (\emph{cf}.\ \cite[Theorem A]{ACV20} or \cite[Theorem 2.1]{ACV19}), $J^p_{a,X_\C}$ admits a \textbf{distinguished model} $\mathscr{J}^p$ over~$k$, that is, an abelian variety $\mathscr{J}^p$ over $k$ such that 
	\begin{equation}\label{distinguished model of J^p_{a,X_C}}
		\mathscr{J}^p\otimes_k \C\cong J^p_{a,X_\C}
	\end{equation} is an isomorphism, and it is \emph{unique} (hence the name ``distinguished'') in that the above Abel--Jacobi map  is $\mathrm{Aut}(\C/k)$-equivariant (see the first diagram in \cite[\S 2.2]{ACV19} for the precise meaning).  The next step is to show that $\mathscr{J}^p$ can be extended to an abelian scheme over the given $S$. This is accomplished in \cite[\S 5]{ACV19}, the resulting object is the abelian scheme in Theorem \ref{ACV, Main Theorems}, and its generic fiber is the above distinguished  model~$\mathscr{J}^p$.
\end{void}

\begin{void}[Distinguished normal function]\label{Generic normal functions} 
	Let $X/S$ be a smooth projective family of relative dimension~$d$ over a smooth $\C$-variety $S$. Note that  every Zariski open $U\subset S$ is an object $U\in \mathsf{Sm}/S$ (Notation \ref{Notation for ACV's functor}(2)). For every cycle $z\in \CH^p_\alg(X_U/U)$ we have an algebraic normal function 
	\[
	v_{z,a}=\Phi^p_{\AJ}(U)(z)\colon U\to J^p_{a,X_U/U}=J^p_{a,X/S}|_U
	\] (where the last identification follows from Theorem \ref{ACV, Main Theorems}(1)). Let $\eta\in S$ be the generic point so that it is contained in every nonempty Zariski open $U\subset S$. Varying such $U$, we obtain a morphism
	\[
	\nu_{z_\eta}\colon \Spec k\to J^p_{a,X_{k}}.
	\] In fact, this is called the \textbf{distinguished normal function} associated to $z\in \CH^p_\alg(X/S)$ in \cite[\S 2.2]{ACV19}, it enjoys certain $\mathrm{Aut}(\C/k)$-equivariance. We also call it the \textbf{generic normal function} (associated to $z$). Passing to the geometric generic fiber, we have
	\begin{equation}\label{generic normal function}
		\nu_{z_{K}} \colon \Spec K\to J^p_{a,X_{K}},
	\end{equation}
	which will be call the  \textbf{geometric generic (algebraic) normal function} (associated to $Z$).  Consequently, we have the \textbf{geometric Abel--Jacobi map}
	\begin{equation}\label{Geometric Abel-Jacobi}
		\AJ^p\colon \CH^p_\alg(X_K)\to J^p_{a,X_K}(K),\quad z\mapsto \nu_z.
	\end{equation} 
\end{void}

\begin{defn}\label{Abel--Jacobi map, alg. trivial}
	Let $X$ be a nice $K$-variety. A cycle class $z\in \CH^p_\alg(X)$ on $X$ is called \textbf{Abel--Jacobi equivalent to zero} or \textbf{Abel--Jacobi trivial}, if $\AJ^p(z)=0$, where  $\AJ^p$ is the geometric Abel--Jacobi map from \eqref{Geometric Abel-Jacobi}.
\end{defn} 
\begin{defn}\label{Incidence equivalence}
	Let $F$ be an algebraically closed field (of any characteristic) and $X$ a nice $F$-variety of dimension~$d$. For $p\in\{1,\ldots,d\}$ consider the subgroup
\begin{align*}
	I^p(X)
	&\coloneq 
	\left\{ z \in \CH^p_{\alg}(X)\;\middle|\;
	\forall\, \text{nice }F\text{-varieties }T,\ 
	\forall\, \alpha \in \CH^{d+1-p}(X \times T): \right. \\
	&\qquad\left. \alpha_* z = 0 \in \CH^1(T) \right\}.
\end{align*} Elements of $I^p(X)$ are called \textbf{incidence equivalent to zero} or \textbf{incidentally trivial} cycles.
\end{defn}
\begin{lem}\label{higher Picard and algebraically trivial cycles}
	Let $F$ be an algebraically closed field and $X$ be a nice $F$-variety of dimension $d$. For every $p\in \{1,\ldots,d\}$, the  universal map $\theta^p$ from Definition \ref{Def. Picard variety} induces an isomorphism
	\[
	\theta^p\colon \CH^p_\alg(X)/I^p(X)\cong \Pic^p_{X/F}(F),
	\] i.e.\ $I^p(X)=\ker\theta^p$.
\end{lem}
\begin{proof}
	This is \cite[Theorem 3.6(i)]{HSaito}.
\end{proof}
\begin{lem}\label{AJ trivial implies I-trivial, complex}
	Let $X$ be a nice $\C$-variety of dimension $d$.  Every Abel--Jacobi trivial cycle on $X$ is incidentally trivial on $X$. In particular, by the previous Lemma \ref{higher Picard and algebraically trivial cycles}, for  every $p\in\{1,\ldots, d\}$ there is a surjective $\C$-morphism of abelian varieties
	\begin{equation}\label{Saito's h^p}
		h^p\colon J^p_{a,X}\to \Pic^p_{X/\C}
	\end{equation}  such that the following diagram commutes:
	\begin{equation}
	\begin{tikzcd}\label{h^pAJ^p=theta^p}
		\CH^p_\alg(X)\ar{r}{\AJ^p} \ar[swap, twoheadrightarrow]{rd}{\theta^p} & J^p_{a,X}(\C) \ar[twoheadrightarrow]{d}{\pi^p}\\
		& \Pic^p_{X/\C}(\C)
	\end{tikzcd}  
\end{equation}
\end{lem}
\begin{proof}
	Let $T$ be a nice $\C$-variety, $\alpha\in \CH^{d+1-p}(X\times T)$ a correspondence. There is natural morphism $[\alpha]\colon J^p_{a,X}\to J^1_T(=J^1_{a,T})$ of complex abelian varieties, making the following diagram
	\[
	\begin{tikzcd}
		\CH^p_\alg(X)\ar{d}{\AJ^p} \ar{r}{\alpha_*} & \CH^{1}_\alg(T)\ar{d}{\cong}\\
		J^p_{a,X}(\C)\ar{r}{[\alpha]} & J^{1}_{T}(\C)
	\end{tikzcd}
	\] commute, see e.g.\ \cite[Theorem 12.24]{LewisBook}, from which the first part of the claim follows. 
\end{proof}
Recall Notation \ref{notation for fields for algebraic intermediate Jacobian}.
\begin{lem}[Chow's rigidity theorem]\label{Chow's rigidity}
	Let $A, B$ be abelian varieties over $K$. Every $\C$-morphism $A_\C\to B_\C$ can be descended uniquely to a $K$-morphism $A\to B$.
\end{lem}
\begin{proof}
	Apply \cite[Theorem 3.19]{Con06} to  the regular field extension $K\subset \C$.
\end{proof}
\begin{rem}[Descent of complex higher Picard varieties]\label{Picard var and base change}
	For $p\in\{1,\ldots, d\}$, $\Pic^p_{X/K}$ is the algebraic representative in codimension $p$ for  \emph{Picard homomorphisms}.   Since $\Pic^p_{X/K}$ and $\Pic^p_{X_\C/\C}$ are algebraic representatives for the Picard homomorphisms  on  $\CH^p_\alg(X)$ resp.\ on $\CH^p_\alg(X_\C)$,  we can apply \cite[Theorem 1(i)]{ACV23} and find
	\[
	\Pic^p_{X_\C/\C}\cong \Pic^p_{X/K}\otimes_K \C.
	\] See especially the proof of \cite[Theorem 5.6]{ACV23} (which uses the notion of $\C/K$-traces of an abelian variety).
\end{rem}
\begin{cor}\label{J^p to Pic^p}
	Let $X$ be a nice $K$-variety  of dimension $d$. For every $p\in \{1,\ldots,d\}$, there is a unique surjective morphism
	\begin{equation}
		h^p\colon J^p_{a,X}\twoheadrightarrow \Pic^p_{X/K}
	\end{equation} such that the base change $h^p_\C$ yields the same morphism from Lemma \ref{AJ trivial implies I-trivial, complex}.
\end{cor}
\begin{proof}
	By formula \eqref{distinguished model of J^p_{a,X_C}}, Lemma \ref{AJ trivial implies I-trivial, complex} and Remark \ref{Picard var and base change}, there is  a $\C$-morphism
	\[
	J^p_{a,X}\otimes_K \C=J^p_{a,X_\C} \to \Pic^p_{X_\C}=\Pic^p_{X_K}\otimes_K \C.
	\] We conclude by Lemma \ref{Chow's rigidity} that there is a $K$-morphism 	$h^p\colon J^p_{a,X}\to \Pic^p_{X/K}$. The surjectivity can be descended (along fpqc morphisms), hence $h^p$  is surjective.
\end{proof}
\begin{void}\label{K-points of Jac and Pic}
	We thus obtain a commutative diagram
	\begin{equation}
		\begin{tikzcd}\label{h^pAJ^p=theta^p}
			\CH^p_\alg(X)\ar{r}{\AJ^p} \ar[swap, twoheadrightarrow]{rd}{\theta^p} & J^p_{a,X_{K}}({K}) \ar[twoheadrightarrow]{d}{h^p(K)}\\
			& \Pic^p_{X}(K).
		\end{tikzcd}  
	\end{equation} For the distinguished normal function $\nu_z$ associated to $z\in \CH^p_\alg(X)$, we obtain a related commutative  diagram
\begin{equation}\label{J^p_a to Pic^p}
\begin{tikzcd}
	\Spec K \ar{r}{\nu_z}  \ar[swap]{rd}{\theta^p_z} & J^p_{a,X} \ar[twoheadrightarrow]{d}{h^p}\\
	& \Pic^p_{X}
\end{tikzcd}
\end{equation}
where  the arrow $\theta^p_z$ corresponds to the $K$-point  $\theta^p(z)\in \Pic^p_{X/K}(K)$.
\end{void}

\begin{rem}\label{relative higher Picard var}
	Let $X/S$ be a smooth projective family of fiber dimension $d$ over a smooth $\C$-variety $S$. It would be an interesting project to study the existence of the ``relative $p$-th Picard variety'' (for $p\in\{2,\ldots,d-1\}$), denoted temporarily by $\Pic^p_{X/S}$, in  relation with the relative algebraic intermediate Jacobians constructed by Achter--Casalaina-Martin--Vial. For example, one may wonder if the morphism $h^p$ from  \eqref{J^p_a to Pic^p} could be realized as the generic part of a (canonically constructed) morphism $J^p_{a,X/S}\to \Pic^p_{X/S}$. 
\end{rem}
Finally, we remark that  that the surjection $h^p$ is conjectured to be an isogeny when the base field is $\C$ (\cite[\S 5]{HSaito}), or more generally over an algebraically closed field $K$ of a complex function field, where one can speak of Abel--Jacobi maps. This property is in fact equivalent to a conjecture of Griffiths on incidentally trivial cycles (see Problem B on p.17 of \cite{GriffithsProblemB}). 
\begin{conj}[Griffiths]\label{Griffiths conjecutre for incidentally trivial cycles}
	Let $X$ be a nice $K$-variety of dimension $d$ and $p\in \{1,\ldots, d\}$. For every $z\in I^p(X)\subset\CH^p_\alg(X)$, there is a positive integer $N$ such that $Nz\in \ker(\AJ^p)$. 
\end{conj}
\begin{lem}\label{hp isog. and Griffiths conj}
	Griffiths's conjecture  for $\CH^p_\alg(X)$ is equivalent to that $h^p\colon J^p_{a,X}\to \Pic^p_X$ is an isogeny.
\end{lem}
\begin{proof}
	The compatibility $\theta^p=h^p\circ \AJ^p$ implies that
	\begin{equation}\label{ker(hp)=AJ(ker of theta-p)}
		\ker(h^p)(K)=\AJ^p(\ker\theta^p)
	\end{equation}
as abelian groups.
	
Suppose $h^p$ is an isogeny, i.e.\ $\ker(h^p)$ is finite.  Then \eqref{ker(hp)=AJ(ker of theta-p)} immediately gives the Griffiths conjecture \ref{Griffiths conjecutre for incidentally trivial cycles}.
	
On the other hand, \eqref{ker(hp)=AJ(ker of theta-p)}	and Griffiths's conjecture  imply that  $\ker(h^p)(K)$ consists solely of torsion points. Since we work over a field of characteristic zero, this can only happen if $\ker(h^p)$ is a finite group scheme.
\end{proof}
\begin{ex}\label{Example: Griffiths conj}
Let $X$ be a $d$-dimension $\C$-variety.
	\begin{enumerate}
		\item When $X$ is an abelian variety, then $h^p$ is an isogeny for every $p$ (\cite[Proposition 5.2]{HSaito}). More generally, this holds when homological triviality and numerical triviality agree up to torsion (\emph{cf}.\ \cite{KleimanFiniteness}).
		
		\item Recall $h^p$ is surjective. When $X$ is a smooth complete intersection, then $$H^{2p-1}(X,\C)=0\quad\text{for } 2p-1<d.$$ In particular, $J^p_{a,X}=0=\Pic^p_X$ and $h^p=0$ is trivially an isogeny  for all such $p$ (\emph{cf}.\ \cite[(6.1)]{HSaito}). 
		
		\item By a hyperplane section argument, if $\dim X=2p-1$ is odd, to prove $h^\ell$ is an isogeny for any $\ell\in{1,\ldots,\dim X}$, it suffices to prove it for $\ell=p$ (\cite[Remark 5.6]{HSaito}).
		
		\item By a theorem of Murre (\cite[Theorem 2.4]{MurreIncidence}), $h^2\colon J^2_{a,X}\to \Pic^2_{X}$ is always an isogeny.
	\end{enumerate}
Since isogenies can be descended along $K\subset \C$ ($K$ algebraically closed), we obtain the similar statements for $h^p$ in the above examples when $X$ is a nice $K$-variety.
\end{ex}
\section{The geometric height pairing à la Kahn}\label{The geometric height pairing à la Kahn}
We review Kahn's framework \cite{Kahn} on the geometric Beilinson--Bloch  height pairing. For practical reasons, we will specialize his general results to the complex function field case. On the other hand, we think it will be useful to consider Chow groups (of  regular models) with $\Q$-coefficients, whereas Kahn did not, but see his definition $\CH^\bullet(X)^{(0)}$, recollected in Remark \ref{Saturated admissible cycles}. 

In this section, $k$  denotes a complex function field. Let $X$ be  a nice $k$-variety of dimension~$d$.

\begin{defn}\label{Model, Kahn}
	By a \textbf{regular model} $f\colon \X \to B$ of $X\to \Spec k$, we mean that $f$ is a proper flat morphism,  $\X$ is  a smooth $\C$-variety (but the structure morphism $f$ need not be smooth) and $B$ is a smooth (but not necessarily proper) $\C$-variety  whose function field $\C(B)$ is isomorphic to $k$ such that the generic fiber of $f$ is isomorphic to $X$; we fix one such isomorphism $\X_k\cong X$. 
\end{defn}
\begin{rem}
	Every nice $k$-variety admits a regular model  by the resolution of singularities due to Hironaka.
\end{rem}
\begin{void}
	Let $f\colon \X\to B$ be a regular model of $X$. Composing the intersection pairing on the smooth $\C$-variety $\X$ with the proper pushforward $f_*$, we obtain a pairing  (for varying $i, r$)
	\begin{equation}\label{f-pairing}
		\langle,\rangle_\X\colon \CH^i(\X)_\Q\times \CH^{d+r-i}(\X)_\Q\to \CH^r(B)_\Q.
	\end{equation}
\end{void}
Let us fix a regular model $(\X,f,B)$ of $X$ in the following. 	Let $j\colon X\hookrightarrow \X$ denote the natural monomorphism.
\begin{notn} 
	For a morphism $Z\to B$, we let $f_Z\colon \X_Z\to Z$ denote the base change morphism. Typically $Z$ will be a (locally closed) subscheme of $B$. 
\end{notn}
\begin{void}
	Let $b\in B^{(1)}$ be a codimension one point, $\iota\colon Z\coloneq \overline{\{b\}}\hookrightarrow B$ the closed immersion. We have a \textbf{cup product} 
	\[
	\CH^i(\X)_\Q\times \CH_j(\X)_\Q\to \CH_{j-i}(\X_Z)_\Q,\quad (\alpha,\beta)\mapsto \Gamma_\iota^!(\beta\times \alpha), 
	\] where $(-)^!$ denotes the refined Gysin pullback (\cite[\S 6.2]{Fulton}). Taking $j=\dim B+i-1$ and composing the cup product with $(f_Z)_*$, we obtain a pairing
	\begin{equation}
		\langle,\rangle_b\colon \CH^i(\X)_\Q\times \CH_j(\X_Z)_\Q\to \CH_{\dim B-1}(Z)_\Q=\CH^0(Z)_\Q=\Q.
	\end{equation}
\end{void}
\begin{void}\label{admissible cycle class}
	We define 
	\begin{align*}
		&	\CH^i(\X)^0_b\\
		=\,{} & \{\alpha\in \CH^i(\X)_\Q\mid j^*\alpha\in \CH^i_{\num}(X)\text{ and }\langle\alpha,\beta\rangle_b=0 \text{ for all }\beta\in \CH_{\dim B+i-1}(\X_Z)_\Q\}
	\end{align*} (the subscript ``num'' stands for numerically trivial cycles).
	And we put
	\[
	\CH^i(\X)^0=\bigcap_{b\in B^{(1)}}\CH^i(\X)^0_b.
	\] Elements of $\CH^i(\X)^0$ will be called \textbf{admissible cycles}.
\end{void}
\begin{rem}
	By \cite[Remark 2.3(b)]{Kahn}, the condition ``$j^*\alpha\in \CH^i_{\mathrm{num}}(X)$'' can be left out without changing the group $\CH^i(\X)^0_b$.
\end{rem}
The basic idea of defining the (geometric) height pairing of cycles  on $X$ is through the intersection pairing of their ``good'' liftings (if they ever exist) on a regular model of $X$. This implies that the intersection pairing is independent of the lifting. An instance is given in the next
\begin{lem}\label{lemma for independence of liftings}
	Let $\alpha\in \CH^i(\X)^0$. If $\beta\in \CH^{d+1-i}(\X)_\Q$ such that $j^*\beta=0$, then $\langle \alpha,\beta\rangle_\X=0\in \CH^1(B)_\Q$ (recall \eqref{f-pairing}).
\end{lem}
\begin{proof}
	See \cite[Proposition 2.8]{Kahn}. The proof continues to work for Chow groups with $\Q$-coefficients.
\end{proof}
\begin{defn}\label{geometric height pairing}
	Recall  $(\X,f,B)$ is a regular model of $X$ and $j\colon X\hookrightarrow \X$ denotes the natural monomorphism. Let 
	\[
	\CH^i(X)^0_f\coloneq \mathrm{im}(j^*\colon \CH^i(\X)^0\to \CH^i(X))\subset \CH^i(X).
	\] be the image. In other words, $\CH^i(X)^0_f$ consists of cycles that admit admissible liftings (with $\Q$-coefficients).
	
	By Lemma \ref{lemma for independence of liftings}, we have a well-defined pairing (for varying $i\geq 1$)
	\begin{equation}
		\langle,\rangle_f\colon \CH^i(X)^0_f\times \CH^{d+1-i}(X)^0_f\to \CH^1(B)_\Q.
	\end{equation} This will be our \textbf{(global) geometric height pairing}.
\end{defn} 
A priori, the geometric height pairing depends on a  regular model $(\X,f,B)$. However, we have the following independence result.
\begin{thm}[Kahn]\label{Kahn's geometric height pairing}
	Let $(\X_1,f_1,B)$ and $(\X_2,f_2,B)$ be two regular models of $X$ (with the same base $B$). Then $$\CH^i(X)^0\coloneq \CH^i(X)_{f_1}^0=\CH^i(X)_{f_2}^0$$ for all $i\geq 1$. Moreover, for any $\alpha\in \CH^i(X)^0, \beta\in \CH^{d+1-i}(X)^0$, we have
	\[
	\langle\alpha,\beta\rangle\coloneq \langle\alpha,\beta\rangle_{f_1}=\langle\alpha,\beta\rangle_{f_2}\in \CH^1(B)_\Q.
	\] 
\end{thm}
\begin{proof}
	This is \cite[Propositions 3.6 \& 3.8]{Kahn}. The proof works for Chow groups with $\Q$-coefficients.
\end{proof}
\begin{rem}\label{Height Pairing and Correspondence}
	The paring $\langle,\rangle$ behaves well under correspondences. More precisely, suppose $X_1, X_2$ are nice $k$-varieties, $\gamma\in \CH^{\dim X_2}(X_1\times X_2)$ a  correspondence. Then we have a canonical isomorphism
	\[
	\langle \gamma^*\alpha,\beta\rangle=\langle \alpha,\gamma_*\beta\rangle \in \Pic(B)_\Q
	\] for $(\alpha,\beta)\in \CH^i(X_2)^0\times\CH^{\dim X_1-i+1}(X_1)^0$.  The proof is contained in \cite[Proposition 3.8]{Kahn}.
\end{rem}
\begin{rem}\label{Saturated admissible cycles}
	In \cite[\S 5]{Kahn}, Kahn introduces 
	\[
	\CH^i(X)^{(0)}\coloneq \{a\in \CH^i(X)\mid \exists n\in \Z\ohne\{0\}: na\in \CH^i(X)^0\},
	\] where the $\CH^i(X)^0$ in the display uses the Chow groups of a regular model with \emph{integral} coefficients. This $\CH^i(X)^{(0)}$ is at least contained in our $\CH^i(X)^0=\CH^i(X)^0_f$ in the sense of Definition \ref{geometric height pairing} (because we uses more general---those with $\Q$-coefficients---admissible cycles to generate the image $\CH^i(X)^0$).
\end{rem}
The next theorem is a geometric analogue of \cite[Lemma 8.1]{Kue96}.
\begin{thm}[Kahn]\label{Analog of Künnemann's Lemma 8.1}
	For every $i\geq 1$, we have $\CH^i_{\alg}(X)\subset \CH^i(X)^0$.
\end{thm}
\begin{proof}
	By the previous remark, it is enough to show $\CH^i_{\alg}(X)\subset \CH^i(X)^{(0)}$, which is  \cite[Theorem 5.6(a)]{Kahn}. Since we will need   the proof later, we reproduce the key ideas.
	
	Fix a regular model $(\X,f,B)$ of $X$. Let $j\colon X\hookrightarrow \X$ be the natural  monomorphism.
	
	\emph{Step 1}. First assume that $d=\dim X=1$ and that there exist sections $\tilde{c}_0,\, \tilde{c}_1\colon B\to \X$ of $f$. Let $c_i$ be the generic fibers of the sections and set $a\coloneq [c_0]-[c_1]\in \CH^1(X)$. 
	
	\emph{Claim 1}. There exists an integer $N>0$ such that $Na\in \CH^1(X)^0$. More precisely, setting 
	\[\tilde{a}\coloneq [\tilde{c}_0(B)]-[\tilde{c}_1(B)]\in \CH^1(\X),\] we claim that there is an integer $N>0$ and $\xi\in \ker j^*$ such that $N\tilde{a}+\xi\in \CH^1(\X)^0$. This implies the previous claim by taking the generic fiber of $\tilde{a}$.
	\begin{proof}[Proof of Claim 1]
		The result is classical if $\dim B=1$ from the theory of fibred surfaces (see e.g.\ \cite[Theorem 9.1.23]{Liu}). Kahn verifies that this continues to hold over a higher dimensional base.
	\end{proof}
	
	\emph{Step 2}. Let $z\in \CH^i_\alg(X)$.  By our definition of algebraic triviality, there exist a nice $k$-curve $C$, a correspondence $\alpha\in \CH^i(C\times_k X)$ and $c_0, c_1\in C(k)$ such that
	\[
	z=\alpha_{c_0}-\alpha_{c_1}=([c_0]-[c_1])^!\alpha.
	\]
	
	\emph{Claim 2}. There exists a closed subset $F\subset B$ with $\mathrm{codim}_B(F)\geq 2$ such that $C\to \Spec k$ admits a regular model $g\colon \mathcal{C}\to B\ohne F$ and the $c_i$ lift to sections $\tilde{c}_i\colon B\ohne F\to \mathcal{C}$ of $g$.
	\begin{proof}[Proof of Claim 2]
		See \cite[Lemma 5.8]{Kahn}.
	\end{proof}
	Note that the function field of $B\ohne F$ is still isomorphic to $k$. So by restriction we get a new regular model of $X$  with the extension properties as in \emph{Claim} 2, which will be  denoted again by $(\X,f,B)$.
	
	Let $z=\alpha_{c_0}-\alpha_{c_1}\in \CH^i_{\alg}(X)$ as in Step 2. By Step 1, there exists an integer $N>0$ such that $N([c_0]-[c_1])\in \CH^1(C)^0$. The  trick is to observe that
	\[
	Nz=N(([c_0]-[c_1])^!\alpha)=(\alpha^\top)^*N([c_0]-[c_1])
	\] (Remark \ref{Weil-Bloch}). Then \cite[Proposition 3.6]{Kahn} implies that the right most term lies  in  $\CH^i(X)^0$. This completes the proof.
\end{proof}
\begin{cor}\label{Geometric height pairing for alg. trivial cycles}
	Let $X$ be a nice $k$-variety of dimension $d$,  $B$ be a smooth $\C$-variety with function field $k$, which is the base of a regular model of $X$. There is a well-defined geometric height pairing 
	\[
	\langle,\rangle\colon \CH^i_{\alg}(X)\times \CH^{d+1-i}_{\alg}(X)\to \CH^1(B)\cong \Pic(B)_\Q
	\] for every $i\geq 1$.
\end{cor}
\begin{rem}
	There is an unconditional homological approach to the geometric height pairing of homologically trivial cycles, originally due to Beilinson (\cite[\S 1]{BeilinsonHP}) and recently revisited and generalized by Rössler--Szamuely \cite{RS22} over a higher dimensional base. However, note that $$\CH^\bullet_{\hom}(X)\subset \CH(X)^0$$ is still a conjecture for general $X$, even over a 1-dimensional base (\cite[Conjecture 5.2]{Kahn}). We will not use the homological construction and refer the interested reader to \cite[\S 2E]{Kahn} and the preprint \cite{Wisson} by T. Wisson  for a detailed comparison. 
\end{rem}
\section{The geometric Néron--Tate height pairing}\label{Geometric NT height pairing}
This section is based on Moret-Bailly's classical works \cite{MB_adm} and \cite[\S III.3]{Ast129}. (The arithmetic counterpart is also discussed in loc.\ cit., but we focus on the function field case.) We assume the knowledge of Néron models covered in \cite{BLR}. We will explain  a line bundle-valued pairing of the geometric Néron--Tate pairing. This is similar, in spirit, to  the pairing in  the previous section.  

The notion of rigidification will play an important role in all the extension results of line bundles.
\begin{defn}\label{rigidified line bundle}
	Let $f\colon X\to S$ be a  morphism of schemes and $i\colon S\to X$ be a section of~$f$. A line bundle $L$ on $X$ is \textbf{rigidified along}  $i$, if there is an isomorphism
	\[
	\phi\colon i^*L\cong \O_S.
	\] If $(L,\phi)$ and $(L',\phi')$ are two rigidified line bundles along $i$, we say that they are \textbf{isomorphic}, if there is an isomorphism $u\colon L\cong L'$ which is compatible with the rigidifications:
	\[
	\phi'\circ i^*u=\phi\colon i^*L\cong \O_S.
	\] We thus obtain a category (groupoid) $\mathsf{PicRig}(X,i)$ with objects consists of rigidified line bundles, and morphisms are isomorphisms between rigidified line bundles.
\end{defn}
\begin{lem}\label{equivalence of cat. of rigidified line bundle}
		Let $f\colon X\to S$ be a smooth  morphism of schemes over   a connected regular noetherian scheme $S$. Assume moreover that $f$ has \emph{connected fibers} over all codimension-one points. Let $i$ be a section of~$f$.  Let $U\subset S$ be a dense open subscheme. The restriction $X_U\subset X$ induces an equivalence of categories
		\[
		\mathsf{PicRig}(X,i)\to \mathsf{PicRig}(X_U,i_U),
		\]  where $i_U\colon U\to X_U$ denotes the section induced by $i$.
\end{lem}
\begin{proof}
When $\dim S=1$, this is \cite[Lemme 2.8.2.1]{MB_adm}. The proof idea can be adapted to the general case.  
\end{proof}
We deduce an important extension result of Poincaré bundles. 
\begin{thm}[Prolongation of the Poincaré biextension bundle]\label{Extension of Poincare bundle}
	Let $S\subset \Sbar$ be a dense open subset of a connected regular noetherian scheme $\Sbar$. Let $A\to S$ be  an abelian scheme, and suppose there are smooth commutative group schemes $G/\Sbar$ resp.\ $G'/\Sbar$ extending $A/S$ resp.\ the dual abelian scheme $A^\vee/S$. Assume both $G$ and $G'$  have \emph{connected} fibers over codimension-one points. Then the Poincaré bundle $\PP$ over $A\times_S A^\vee$, rigidified along the identity section,  admits a \emph{unique} extension  $\overline{\PP}$ over $G\times_{\Sbar} G'$, rigidified along the identity section.
\end{thm}
\begin{rem}
	For a product $G\times_S G'$ of commutative group schemes over a scheme $S$, we let 
	\[
	[-1]\colon G\times_S G'\to G \times_S G',\quad (x,y)\mapsto (-x,-y)
	\] denote the total multiplication by $-1$. Now use the notation from Theorem \ref{Extension of Poincare bundle}. The Poincaré  bundle $\PP$ is symmetric, i.e. $[-1]^*\PP\cong \PP$. By Lemma  \ref{equivalence of cat. of rigidified line bundle},  the extension $\overline{\PP}$ is also symmetric.
\end{rem}
Now we specialize  to the case where $S$ is a smooth $\C$-curve; let $k$ be the corresponding  complex function field. We fix an algebraic closure $k\subset K$.
\begin{notn}\label{Notation for Neron model}
	For an abelian variety $A$ over  $k\cong \C(S)$, we let   $\NN(A)\to {S}$ denote the Néron model of $A$ over $S$ and $\NN(A)^0\to S$ be  the identity component. Recall $\NN(A)^0\to S$ is a smooth commutative group scheme of finite type over $S$.
	
	Let $x\colon \Spec K\to A$ be a section of the abelian variety $A$, we write $\tilde{x}\colon {S}\to \NN(A)$ for the unique morphism that extends  $x$, granted by the Néron mapping property.
\end{notn}
\begin{rem}\label{mapping into N^0}
	The component groups of  $\NN(A)$ at each  place $v$ of $\C(S)$ are trivial for all but finitely many $v$. Hence, there is a positive integer $n_x$ such that $n_x\tilde{x}\in \NN^0(A)(S)$.
\end{rem}

\begin{void}
	Let $A$ be an abelian variety over $k$. Let $\PP$ be the  Poincaré bundle over $A\times_k A^\vee$. By Theorem \ref{Extension of Poincare bundle}, $\PP$ extends  uniquely to a  rigidified line bundle $\overline{\PP}$ over $\NN(A)^0\times_{{S}} \NN(A^\vee)^0$.
	
	Let $x\in A(k),\, y\in A^\vee(k)$. By Remark \ref{mapping into N^0}, there is a positive integer $n$ such that $(x, y)\in A(k)\times A^\vee(k)$ extends to a section
	\[
	(n\tilde{x}, n\tilde{y})\colon S\to \NN(A)^0\times_{S} \NN(A^\vee)^0
	\] of $(\NN(A)^0\times_{S} \NN(A^\vee)^0)/S$.
\end{void}
\begin{defn}
	Keep the previous notations. The \textbf{line bundle for the Néron--Tate height}  associated to $(x, y)\in A(k)\times A^\vee (k)$  is the $\Q$-line bundle
	\begin{equation}
		\mathcal{B}_{<x,y>}^k\coloneq \frac{(n\tilde{x}, n\tilde{y})^*{\overline{\PP}}}{n^2}\in \Pic({S})_\Q.
	\end{equation} As a $\Q$-line bundle, 	$\BB^k_{<x,y>}$ is clearly independent of $n$.
	
	If, in addition, $S$ is proper, the \textbf{geometric Néron--Tate height pairing} of $(x, y)\in A(k)\times A^\vee (k)$ is  defined to be the degree of the line bundle
	\begin{equation}\label{geometric NT height pairing}
		\langle x,y\rangle_A\coloneq \deg\mathcal{B}_{<x,y>}^k\in \Q.
	\end{equation} 
\end{defn} 
\begin{rem}[Base change and height]\label{Base change and height}
	Let  $A$ be an abelian variety over $k$. If $x\in A(k), y\in A^\vee(k)$ are geometric points, there is a finite extension $k'/k$ such that $x\in A(k'), y\in A^\vee(k')$. Let $\nu'\colon S'\to S$ be the normalization of $S$ in $k'$ so that $S'$ is also a smooth $\C$-curve. We can form the line bundle $\BB_{<x,y>}^{k'}\in \Pic(S')_\Q$ by the previous definition. We set	 (recall $\nu'$ is a finite morphism)
	\begin{equation}\label{Base change and height, formula}
		\BB_{<x,y>}^{k',k}\coloneq \frac{1}{[k':k]}\nu'_*\BB_{<x,y>}^{k'}\in \Pic(S)_\Q,
	\end{equation} where we identify $\Pic(S')_\Q\cong \CH^1(S')_\Q$ using the regularity of $S'$ and apply the pushforward $\nu'_*$ to the cycle class corresponding to $\BB_{<x,y>}^{k'}$ (and take $\BB_{<x,y>}^{k',k}$ to be the line bundle corresponding to the pushed forward cycle class).

	We claim that the definition does not depend on the choice of $k'$. Indeed, if $k''/k$ is another finite extension such that $x\in A(k''), y\in A^\vee(k'')$, we can form the line bundle $\BB^{k'',k}_{<x,y>}\in \Pic(S)_\Q$ in the same way. Now let $\tilde{S}\to S'\times_S S''$ be the normalization. We have a commutative diagram (the meaning of the notations explains itself)
	\[
	\begin{tikzcd}
		\tilde{S}\ar{d}{\alpha''} \ar{r}{\alpha'} & S' \ar{d}{\nu'}\\
		S''\ar{r}{\nu''}  & S
	\end{tikzcd}
	\] Note that all the arrows are finite morphisms whose degree are given by the corresponding function field extensions. Write $\tilde{k}=\C(\tilde{S})$. By construction, there are natural morphisms $\Spec \tilde{k}\to \Spec k''$ and $\Spec \tilde{k}\to \Spec k'$ such that $x, y$ can be pulled back to the same $\tilde{x}, \tilde{y}\in A(\tilde{k})$ along either of the natural maps. Again we can form the line bundle $\BB^{\tilde{k}}_{<\tilde{x}, \tilde{y}>}\in \Pic(\tilde{S})_\Q$. The commutativity of the above diagram and \cite[\S III, Proposition 3.3.4(i)]{Ast129} imply that
	\[
	\BB^{k'',k}_{<x,y>}\cong \frac{1}{\deg\alpha''\deg \nu''} (\nu''\circ\alpha'')_*\BB^{\tilde{k}}_{<\tilde{x}, \tilde{y}>}\cong\frac{1}{\deg\alpha'\deg\nu'} (\nu'\circ \alpha')_*\BB^{\tilde{k}}_{<\tilde{x}, \tilde{y}>}\cong \BB^{k',k}_{<x,y>}
	\] as $\Q$-line bundles over $S$.
\end{rem}
\begin{notn}\label{Notation of line bundle for geometric NT}
	Continuing the notations of the previous remark, we will denote the line bundle $\BB_{<x,y>}^{k',k}$ henceforth by $$\BB_{<x,y>}\in \Pic(S)_\Q$$ for geometric points $(x,y)\in A(K)\times A^\vee(K)$. 
\end{notn}
\begin{void}
	If $S$ is a nice $\C$-curve, on taking the degree we thus obtain  the   bilinear map of the geometric Néron--Tate height pairing
	\[
	\langle,\rangle_A\colon A(k)\times A^\vee(k)\to \Q.
	\] The construction agrees with the (more common) geometric Néron--Tate height defined using Tate's limit arguments, \emph{cf}.\ \cite[\S 5]{MB_adm}. 
\end{void}

\begin{void}\label{NT pairing and dual}
	Let  $f\colon A\to B$ be a $k$-morphism of abelian varieties. Let  $f^\vee\colon B^\vee\to A^\vee$ be the dual morphism. Let $x\in A(k)$ and $y\in B^\vee(k)$.  It is well-known that the Néron--Tate height pairing behaves well under dual morphisms:
	\begin{equation}\label{NT pairing and dual_num}
		\langle f(x),y\rangle_B=\langle x,f^\vee(y) \rangle_A.
	\end{equation}
	We want to deduce this numerical equality from an isomorphism of line bundles. A similar stance will be adopted throughout this and next sections.
	
	Let $P_A$ resp.\ $P_B$ denote the Poincaré bundle associated to $A$ resp.\ $B$. By the  characterization of  dual morphisms, we have an isomorphism of line bundles over $A\times_k B^\vee$:
	\begin{equation}\label{Poincare bundle and dual morphism, 1}
		(f\times \id_{B^\vee})^*P_B\cong 	(\id_{A}\times f^\vee)^*P_A.
	\end{equation}
	Let us assume that the abelian varieties $A/k$ and $B/k$ spread out to  abelian schemes $\A/S$ and $\mathscr{B}/S$ over a smooth $\C$-curve $S$, let $S\subset \Sbar$ be the smooth projective compactification and let 
	$\NN_A,\, \NN'_A, \, \NN_B,\, \NN'_B$  be the \emph{identity components} of the Néron models of   $\A,\, \A^\vee,\, \mathscr{B},\, \mathscr{B}^\vee$ respectively over $\Sbar$. Let $\PP_A,\, \PP_B$ denote the Poincaré bundles associated to $\A,\, \mathscr{B}$ and $\overline{\PP_A},\ \overline{\PP_B}$ be the extensions  provided by Theorem \ref{Extension of Poincare bundle}.
	
	Without loss of generality, we can spread out the isomorphism \eqref{Poincare bundle and dual morphism, 1} to an isomorphism 
	\begin{equation}\label{Poincare bundle and dual morphism, 2}
		(f\times 1)^*\PP_B\cong  (1\times f^\vee)^*\PP_A
	\end{equation} of rigidified line bundles on the abelian scheme $\A\times_S \mathscr{B}^\vee$, where we abusively use  $f, f^\vee$ to denote the  spread out morphisms and the 1's denote the intended identities. By the Néron mapping property, $f$ and $f^\vee$ extend further to morphisms
	\[
	\tilde{f}\colon \NN_A\to \NN_B,\quad \tilde{f}^\vee\colon \NN'_{B}\to \NN'_A.
	\] Note that $(\tilde{f}\times 1)^*\overline{\PP_B}$ resp.\ $(1\times \tilde{f}^\vee)^*\overline{\PP_A}$ extend 	$(f\times 1)^*\PP_B$ resp.\ $(1\times f^\vee)^*\PP_A$ in the category of \emph{rigidified} line bundles. Since extensions of rigidified line bundles are unique, \eqref{Poincare bundle and dual morphism, 2} induces an isomorphism
	\begin{equation}\label{Poincare bundle and dual morphism, 3}
		(\tilde{f}\times 1)^*\overline{\PP_B} \cong (1\times \tilde{f}^\vee)^*\overline{\PP_A},
	\end{equation} of rigidified line bundles over $\NN_A\times_{\Sbar} \NN'_B$.
	
	Finally, let 
	\[
	(\tilde{x},\tilde{y})\colon \Sbar\to \NN_A\times_{\Sbar} \NN'_B
	\] be  the extended sections, where to simplify notation we have suppressed the possibly necessary positive multiples in front of the extended sections. Pulling back the isomorphism \eqref{Poincare bundle and dual morphism, 3} further along $(\tilde{x},\tilde{y})$, we find 
	\begin{equation}
		\left(\widetilde{f(x)},\tilde{y}\right)^*\overline{\PP_B}\cong \left(\tilde{x},\widetilde{f^\vee(y)}\right)^*\overline{\PP_A},
	\end{equation}
	noting that $\tilde{f}\circ \tilde{x}=\widetilde{f(x)}$ and  $\tilde{f}^\vee\circ \tilde{y}=\widetilde{f^\vee(y)}$ by $B(k)=\NN_B(\Sbar)$ and $A^\vee(k)=\NN'_A(\Sbar)$ (recall Notation \ref{Notation for Neron model}). On taking the degree in the last display, we recover \eqref{NT pairing and dual_num}. 
\end{void}
\begin{rem}\label{field of definition of morphism, Poincare bundle}
	More generally, suppose  we have abelian varieties $A, B$ defined over an algebraic closure $K$ and suppose there is a $K$-morphism $A_{K}\to B_{K}$ of the base changes. Consider geometric points $x\in A(K), y\in B^\vee(K)$. Since the given $K$-morphism is of finite type (it is even projective), there is a finite extension $k'/k$ such that
	\begin{itemize}
		\item  $A_{K}\to B_{K}$ descends to a $k'$-morphism $f'\colon A_{k'}\to B_{k'}$ and
		\item $x\in A(k'), y\in B^\vee(k')$.
	\end{itemize} Let $\Sbar$ be a smooth $\C$-curve such that $\C(\Sbar)\cong k$. Let  $\Sbar'\to \Sbar$ be the normalization of $\Sbar$ in~$k'$. Using a similar argument as in \ref{NT pairing and dual}, we find
	\begin{equation}
		\left(\widetilde{f'(x)},\tilde{y}\right)^*\overline{\PP'_B}\cong \left(\tilde{x},\widetilde{f'^\vee(y)}\right)^*\overline{\PP'_A} \in \Pic(\Sbar')_\Q
	\end{equation} with the suggestive notation. We then use Remark \ref{Base change and height} to push the above isomorphism down to an isomorphism in $\Pic(\Sbar)_\Q$. 
\end{rem}
\section{Comparison of the line bundles}\label{Comparison of the (extended) Poincaré bundles}
In this section we prove our main results.
\begin{void}
	Let $k$ be the function field of a smooth $\C$-variety $S$ (here $S$ can be of arbitrary dimension).  We fix a regular model (Definition \ref{Model, Kahn}) $\X\to S$ of a nice variety $X\to \Spec k$. By Corollary \ref{Geometric height pairing for alg. trivial cycles}, there is a well-defined pairing
	\[
	\CH^p_\alg(X)\times \CH^q_\alg(X)\to \Pic({S})_{\Q},
	\] where $p+q=\dim X+1$. Let us denote the  line bundle of  this pairing  by $\E_{<x,y>}$ for $(x,y)\in \CH^p_\alg(X)\times\CH_\alg^q(X)$. We call it the  \textbf{line bundle for the geometric height pairing} associated to $(x,y)$.

	If $S$ is  nice $\C$-curve, we can speak of the \textbf{geometric Beilinson--Bloch height pairing} of $x,\, y$ given by
	\[
	\langle x,y\rangle_X\coloneq \deg\E_{<x,y>}\in \Q.
	\] We thus obtain a bilinear map of the (global) geometric Beilinson--Bloch height
	\[
	\langle,\rangle_X\colon	\CH^p_\alg(X)\times\CH^q_\alg(X)\to \Q
	\] for algebraically trivial cycles.
\end{void}
\begin{rem}\label{Base change and height, Kahn}
Parallel to Remark \ref{Base change and height} and using a similar argument, we can also define the $\Q$-line bundle $$\E_{<x,y>}\in \Pic(S)_\Q$$ for $(x,y)\in \CH_\alg^p(X_{K})\times \CH^q_\alg(X_{K})$.
\end{rem}
\begin{ex}
	Let $C$ be a nice curve over a complex function field $k$ of a smooth $\C$-curve~$S$. Let $x,y\in \CH^1_\alg(C)$, which are just degree-zero divisor classes. Let $\pi\colon\mathcal{C}\to \Sbar$ be a semistable model of $C$, in particular, this is a regular model of $C\to \Spec k$ in our sense (Definition \ref{Model, Kahn}). By the theory of admissible pairing (\emph{cf}. S.W. Zhang's paper \cite{ZhangAP}), the cycle classes $x, y$ admit admissible extensions $x^h, y^h\in \CH^1(\mathcal{C})^0\subset \CH^1(\mathcal{C})_\Q$ (notations as in \ref{admissible cycle class}), and we have an isomorphism of $\Q$-line bundles over $\overline{S}$:
	\[
	\E_{<x,y>}\cong \pi_*\langle x^h,y^h\rangle,
	\] where the right hand side denotes the Deligne pairing $\Pic(\mathcal{C})_\Q\times\Pic(\mathcal{C})_\Q\to \Pic(\overline{S})_\Q$.
\end{ex}
\begin{notn}
	Let $C$ be a nice curve over  $k$ with the Jacobian $J=\Pic^0_{C/k}$. and 
	\[
	\AJ\colon \CH^1_\alg(C_{K})\to J(K)
	\] be the (classical) Abel--Jacobi map (which is an isomorphism). Let $\lambda_C^1\colon J\to J^\vee$
	be the principal polarization from Example \ref{Polarization on curves}; it differs from the usual one using the theta divisor by a minus sign. Note that $\AJ$ is definable over $k\subset K$.
\end{notn}
The Faltings--Hriljac formula (see e.g.\ \cite[\S 6]{Kue96}) addresses the natural question about  relations between the Beilinson--Bloch height pairing $\langle x,y\rangle_C$  and the Néron--Tate height pairing $\langle \AJ(x),\lambda^1_C(\AJ(y))\rangle_J$ for $x,y\in \CH^1_\alg(C_{K})$. In the current geometric setup, we have a refined   result observed by Holmes--de Jong.
\begin{thm}[Geometric Faltings--Hriljac formula]\label{Faltings-Hriljac}
	Let $C$ be a nice curve over a complex function field $k$ of a \emph{nice} curve $S$. Let $x,y\in \CH^1_\alg(C_{K})$. There is an isomorphism of $\Q$-line bundles over ${S}$
	\[
	\E_{<x,y>}\cong \BB_{<\AJ(x),\lambda^1_C(\AJ(y))>}.
	\] In particular, by taking the degree, we obtain
	\[
	\langle x,y\rangle_C=\langle\AJ(x),\lambda^1_C(\AJ(y)) \rangle_J\in \Q.
	\]
\end{thm}
\begin{proof}
	This follows from \cite[Proposition 4.3]{HdJ15}. We note that the proof uses the description of the Deligne  pairing through the determinant of cohomology.
\end{proof}
\begin{rem}
	Keep the notations of the above theorem. \cite[Theorem 5.1]{HdJ15}, which is a special case of Theorem \ref{continuous extension of the Poincaré bundle}, implies that  $\BB_{<\AJ(x),\lambda^1_C(\AJ(y))>}$ is further isomorphic to the Lear extension of the biextension bundle $\BB_{x,y}=(\nu_x,\nu_y^\vee)^*\PP$. The isomorphism of these three   line bundles is the prototype of our main result Theorem \ref{Main Result, alg. trivial cycles}.
\end{rem}
\subsection{Geometric height vs.\ Néron--Tate height}
\begin{void}
	Let $X$ be a nice $k$-variety over a complex function field $k$ of a smooth curve $S$. Recall the positive integer $k^p_X\coloneq k^p_{X_{K}}$ from Lemma \ref{k^p_X}. Recall the isogeny $$\lambda^q\coloneq \lambda^q_{X_{K}}\colon \Pic^q_{X_{K}}\to \Pic^{p,\vee}_{X_{K}}$$ appearing in \eqref{Pic(alpha^*)}, where $p+q=\dim X+1$ as before. Recall also the universal map $$\theta^p\colon \CH^p(X_{K})\to \Pic^p_{X_{K}}(K)$$ to the $p$-th Picard variety. 
\end{void}  
\begin{rem}\label{theta^p and field of definition}
	Let $u\in \CH^p_\alg(X)$ be an algebraically trivial cycle class over $k$. Recall we have only defined  $\theta^p$ for cycles  defined over an algebraically closed field. Formally, we can only apply $\theta^p$ to the base change $\bar{u}\in \CH^p_\alg(X_{K})$ over a chosen algebraic closure $k\subset K$. However, note that the $p$-th Picard variety $\Pic^p_{X_{K}/K}$ descends to an abelian variety $\Pic^p_{X/k}$ over $k$ (Remark \ref{Higher Picard var over perfect field}). We write $\theta^p(u)\in \Pic^p_{X/k}({K})$ for $\theta^p(\bar{u})$ under the identification $\Pic^p_{X_{K}/K}(K)=\Pic^p_{X/k}(K)$. Using Remark \ref{Base change and height}, we can form the $\Q$-line bundle of the geometric Néron--Tate height pairing on $\Pic^p_{X/k}$. A similar practice will also be applied to the geometric Beilinson--Bloch height pairing, justified by Remark \ref{Base change and height, Kahn}.
\end{rem}
With these preparations,  we can now run the proof of the said analogue of \cite[Theorem 8.2]{Kue96}.
\begin{thm}\label{Kue96, Thm.8.2}
	Let $\X\to S$ be a regular model (Definition \ref{Model, Kahn}) of a nice $k$-variety $X$. Let $z\in \CH^p_\alg(X), w\in \CH^q_\alg(X)$ be algebraically trivial cycles over $k$ with $p+q=\dim X+1$.  There is an isomorphism 
	\[
	\E_{<z,w>}\cong (k^p_X)^{-1}\BB_{<\theta^p(z),\lambda^q\circ \theta^q(w)>}
	\]
	of $\Q$-line bundles over the  smooth $\C$-curve ${S}$.
\end{thm}
\begin{proof}
	By the  proof of Theorem \ref{Analog of Künnemann's Lemma 8.1}, we can find a nice $k$-curve $C$,  $\gamma\coloneq [c_0]-[c_1]\in \CH^1_\alg(C)$ with $c_0, c_1\in C(k)$ and a correspondence $\alpha\in \CH^p(C\times_{k} X)$ such that $z=\alpha_*\gamma$ as cycle classes on $X$.
	
	We have a chain of isomorphisms of line bundles over $S$:
	\begin{align*}
		\E_{<z,w>}&\cong \E_{<\alpha_*\gamma,w>}\\ 
		&\cong \E_{<\gamma,\alpha^*w>} & (\text{by  Remark \ref{Height Pairing and Correspondence}})\\ 
		&\cong \BB_{<\AJ(\gamma),\lambda^1_C\circ\AJ(\alpha^*w)>} & (\text{by Theorem \ref{Faltings-Hriljac}}) \\  
		&\cong  \BB_{<\AJ(\gamma),\lambda^1_C\circ\Pic(\alpha^\top)\circ\theta^q(w)>} & (\text{by  } \eqref{functoriality of Pic()})\\
		&\cong (k^p_{X})^{-1}\BB_{<\AJ(\gamma),\Pic(\alpha)^\vee\circ \lambda^q\circ \theta^q(w)>} & (\text{by Lemma \ref{properties of Poincare cycles}(3)})\\
		&\cong (k^p_{X})^{-1}\BB_{<\Pic(\alpha)\circ \AJ(\gamma),\lambda^q\circ \theta^q(w)>} & (\text{by  \eqref{NT pairing and dual}})\\
		&\cong  (k^p_{X})^{-1}\BB_{<\theta^p(z),\lambda^q\circ \theta^q(w)>} &  (\text{by  } \eqref{functoriality of Pic()}),
	\end{align*} where in each step we  use a similar argument as conducted in  \ref{NT pairing and dual}.
\end{proof}
\subsection{Néron--Tate height vs.\ biextension height}
\begin{void}[Dual tori]
	Let $X$ be a compact complex manifold of dimension $d$. For every $p\in\{0,\ldots,d\}$ we have the Poincaré duality
	\[
	H^{2p}(X,\C)\cong H^{2d-2p}(X,\C)^\vee.
	\] Now let $X\to S$ be a projective submersive (holomorphic) map of complex manifolds of relative dimension $d$. For every  $p\in \{1,\ldots,d\}$ the Poincaré duality induces a natural isomorphism
	\begin{equation}\label{Dual tori}
		\PD\colon J^p_{X/S}\cong J^{d+1-p,\vee}_{X/S}
	\end{equation} of relative tori in the holomorphic category.
\end{void}
\begin{void}[Algebraic Poincaré dual]
	Now let $X\to S$ be a smooth projective morphism. By Theorem \ref{ACV, Main Theorems}, we have a relative algebraic subtorus $(J^p_{X/S})^\an\subset J^p_{X^\an/S^\an}$. Let us simplify the notations by ignoring the $(-)^\an$ everywhere. By definition of dual tori, the inclusion $J^p_{a,X/S}\subset J^p_{X/S}$ induces a natural surjective morphism 
	\begin{equation}\label{dual morphism of intermediate Jacobians}
		J^{p,\vee}_{X/S}\twoheadrightarrow J^{p,\vee}_{a,X/S}
	\end{equation} ($J^{p,\vee}_{a,X/S}$ denotes the dual abelian scheme). Set $q=d+1-p$. Using the Poincaré duality, we have a holomorphic map
	\begin{equation}\label{PD_a}
		J^q_{a,X/S}\hookrightarrow J^q_{X/S}\cong J^{p,\vee}_ {X/S}\twoheadrightarrow J^{p,\vee}_{a,X/S}.
	\end{equation}  Since the domain and the codomain are projective, by the GAGA principles, the map \eqref{PD_a} is indeed algebraic. On restricting to the generic fibers, \eqref{PD_a} induces a natural (algebraic) morphism 
	\begin{equation}\label{J^p_a,X to J^q,dual_a,X}
		\PD_a\colon 	J^p_{a,X_k}\to J^{q,\vee}_{a,X_k}
	\end{equation} of abelian varieties over the complex function field $k$. 
\end{void}

	Let $X$ be a nice variety over a complex function field $k$ of a smooth $\C$-curve $S$; let $k\subset K$ be an algebraic closure. So far we have introduced two types of abelian varieties: the higher Picard variety $\Pic^p_{X_K}$ and the algebraic intermediate Jacobian $J^p_{a,X_K}$. As before we choose $p,q\in \N$ such that $p+q=\dim X+1$.

\begin{void}
	We shall write $\nu_{\bar{z}}$ for the geometric generic normal functions  associated to $z\in \CH^p_\alg(X)$ (Remark \ref{Generic normal functions}). Recall the isogeny $\lambda^q_X\colon \Pic^q_{X_{K}} \to \Pic^{p,\vee}_{X_{K}}$ induced by a Poincaré $p$-cycle (\emph{cf}. \eqref{Pic(alpha^*)}). Define 
	\[
	\nu_{\bar{w}}^\vee\colon \Spec K\to J^q_{a,X_{K}}\xrightarrow{\PD_a} J^{p,\vee}_{a,X_{K}}
	\] for $w\in \CH^q_\alg(X)$. 
\end{void}
\begin{notn}
	Let us denote by
	\[
	P_a\to \Pic^p_{X}\times_{K} \Pic^{p,\vee}_{X}\quad\text{resp.}\quad P_h\to J^p_{a,X}\times_{K} J^{p,\vee}_{a,X}
	\]  the respective Poincaré line bundles.   Recall $\Pic^p_{X_{K}}\cong\Pic^p_{X}\otimes_K K$ and $J^p_{a,X_{K}}\cong J^p_{a,X}\otimes_K K$ are compatible with base change w.r.t.\ the  algebraic closure $k\subset K$.
\end{notn}
The goal of this subsection is to compare these two line bundles.

\begin{void}
	We established   a natural morphism $h^p\colon J^p_{a,X_{K}}\to \Pic^p_{X_{K}}$ in Corollary \ref{J^p to Pic^p}. Let $h^{p,\vee}\colon \Pic^{p,\vee}_{X_{K}}\to J^{p,\vee}_{a,{X_{K}}}$ be the dual morphism (defined over $K$). Let $z\in \CH^p_\alg(X), w\in \CH^q_\alg(X)$. Consider the diagram (where 1 stands for the intended identity map)
	\begin{equation}\label{normal functions and theta}
		\begin{tikzcd}
			& & & \Pic^p_{X_K}\times \Pic^{p,\vee}_{X_K} \\
			\Spec K
			\ar{r}{(\nu_{\bar z},\theta^q_{\bar w})}
			\ar[swap, bend left=15mm]{urrr}{(\theta^p_{\bar z},\lambda_X^q\circ \theta^q_{\bar w})}
			\ar[bend right=15mm]{drrr}{(\nu_{\bar z},k^p_X\nu_{\bar w}^\vee)}
			& J^p_{a,X_K}\times \Pic^q_{X_K}
			\ar{r}{1\times\lambda^q}
			& J^p_{a,X_K}\times \Pic^{p,\vee}_{X_K}
			\ar{ur}{h^p\times 1}
			\ar[swap]{dr}{1\times h^{p,\vee}}
			& \\
			& & & J^p_{a,X_K}\times J^{p,\vee}_{a,X_K}
		\end{tikzcd}
	\end{equation}
Here is where  Griffiths's Conjecture \ref{Griffiths conjecutre for incidentally trivial cycles} enters---We claim that the diagram is commutative if the conjecture holds. By Corollary \ref{J^p to Pic^p}, this is clear for the upper part (without the conjecture). So we are left to argue the lower part. The proof is technical and will be carried out in Section \ref{Section: G_m-biextensions}.
\end{void}

To continue the comparison, we will need an auxiliary extension result:
\begin{thm}\label{continuous extension of the Poincaré bundle}
	Let $A\to S$ be an abelian scheme over a smooth $\C$-curve $S$. Let $S\subset \Sbar$ be a smooth partial compactification. Then the smoothly metrized Poincaré bundle $(\PP,||.||)$  extends  to a continuously metrized $\Q$-line bundle $(\overline{\PP},||.||)$ over  $\NN(A)^0\times_{\Sbar} \NN(A^\vee)^0$.
\end{thm}
\begin{proof}[Comment]
The basic idea is the following. Let $\overline{\PP}$ be the algebraic extension of $\PP$ over  $\NN(A)^0\times_{\Sbar} \NN(A^\vee)^0$ (Theorem \ref{Extension of Poincare bundle}). One can verify that a similar construction in \cite[Theorem 6.1.1]{YZ21}  can be adapted to give $\overline{\PP}$ a structure  of an archimedean adelic line bundle (\cite[\S 3.2]{BurgosKramer}), although $\NN(A)^0\to \Sbar$ is only quasi-projective in general. By definition or construction, such an adelic line bundle comes with a continuous metric (being a uniform limit  of continuous model metrics  over compact subsets). Details will appear in \cite{CZL26}.
\end{proof}
\begin{void}\label{spread out normal function}
	Let $X$ be a nice variety over a complex function field $k$ of a smooth curve $S$,  $z\in \CH^p_\alg(X), w\in \CH^q_\alg(X)$. Let $\X\to S$ be a smooth projective morphism with generic fiber $X\to \Spec k$.  Let $\A\coloneq J^p_{a,\X/S}$ be the relative algebraic intermediate Jacobian over~$S$.   Let $\PP_\A$ be the Poincaré bundle of $\A$. Recall the biextension bundle $$\BB_{z,w}=(\nu_z,\nu_w^\vee)^*\PP_\A,$$ where we spread out the generic normal functions $\nu_z, \nu_w$ to  sections of $J^\bullet_{a,X/S}/S$ (after possibly shrinking $S$) without changing the notation. Note that $\PP_\A$ carries a smooth canonical metric, which induces the biextension metric $||.||$ on $\BB_{z,w}$. Let $S\subset \Sbar$ be a smooth partial compactification. By Theorem \ref{continuous extension of the Poincaré bundle}, the Poincaré bundle $\PP_\A$ extends  as a continuously metrized line bundle $\overline{\PP_\A}$ over  $\NN(\A)^0\times_{\Sbar} \NN(A^\vee)^0$. Furthermore, by Remark \ref{mapping into N^0}, there is a positive integer $n$ such that $\nu_z, \nu_w^\vee$ extend to a morphism 
	\[
	(n\widetilde{\nu_z}, n\widetilde{\nu_w^\vee})\colon \Sbar\to \NN(\A)^0\times_{\Sbar} \NN(\A^\vee)^0.
	\] 
\end{void}
\begin{cor}\label{Lear extension and Neron model}
	Keep the previous notations. Recall Notation \ref{Notation for Neron model}.  The Lear extension of $(\BB_{z,w},||.||)$ is given by
	\[
	[\BB_{z,w},||.||]_{\Sbar}\cong \frac{\left(n\widetilde{\nu_z},n\widetilde{\nu_w^\vee}\right)^*\overline{\PP_\A}}{n^2}\stackrel{\mathrm{def}}{=}\BB_{<\nu_z,\nu_w^\vee>}\in \Pic(\Sbar)_\Q.
	\] 
\end{cor}
\begin{proof}
	The right hand side is a continuously metrized line bundle. The isomorphism then follows from the uniqueness of the Lear extension. 
\end{proof}
\begin{thm}\label{comparison biextension vs NT-height line bundles}
Assume Conjecture \ref{Griffiths conjecutre for incidentally trivial cycles}.	We have an isomorphism of $\Q$-line bundles over  $\Sbar$:
	\[
	k^p_X[\BB_{z,w},||.||]_{\Sbar}\cong \BB_{<\theta^p(z),\lambda_X^q\circ \theta^q(w)>}.
	\] Recall the right hand side is defined in Notation \eqref{Notation of line bundle for geometric NT}.
\end{thm}
\begin{proof}
	Again we write $\bar{u}\in \CH^\bullet_\alg(X_{K})$ for the base change of $u\in \CH^\bullet_\alg(X)$. By the commutativity of \eqref{normal functions and theta}, we find 
	\begin{equation}\label{desired commutativity for comparison of Poincare bundles}
		h^p\circ \nu_{\bar{z}}=\theta^p_{\bar{z}}\in \Pic^p_{a,X}(K),\quad k^p_X\nu_{\bar{w}}^\vee=h^{p,\vee}\circ \lambda^q_X\circ \theta^q_{\bar{w}}\in J^{p,\vee}_{a,X}(K)
	\end{equation} as geometric points. There is a finite extension $k'/k$ such that  the above equalities holds for the base changes $z', w'\in \CH^\bullet_\alg(X_{k'})$ of $z, w$. After possibly shrinking $S\subset \Sbar$ we may assume that 
	\begin{itemize}
		\item $J^p_{a,X}$ resp.\ $\Pic^p_{X}$   spread out to  abelian schemes $\A$ resp.\ $\mathscr{B}$ over $S$;
		\item the generic normal functions $\nu_z, \nu_w$ (defined over $k$) extend to a section of $J^\bullet_{a,\X/S}$ over $S$, where $\X/S$ is a smooth projective morphism with generic fiber $X$.
	\end{itemize}  
	
	Let $\nu\colon \overline{S}'\to \Sbar$ be the normalization of $\overline{S}$ in $k'$ and consider the open subset $S'\coloneq \Sbar'\times_{\Sbar} S$ of $\Sbar'$. Write $\X'\coloneq \X\times_S S'$ and $\A'\coloneq\A_{S'}\coloneq J^p_{a,\X}\times_{S} S'=J^p_{a,\X'/S'}$ (Theorem \ref{ACV, Main Theorems}(1)); the latter admits a finite morphism to $\A=J^p_{a,\X/S}$.
	Consider the base change  of $(\nu_z,\nu_w^\vee)$ over~$S'$:
	\[
	(\nu'_z,\nu'^\vee_{w})\colon S'\to \A'\times_{S'}\A'^\vee.
	\] This gives a {pullback} diagram
	\[\label{pullback, P_h and P'_h}
		\begin{tikzcd}
			S'\ar[swap]{d}{\nu|_{S'}}\ar{r}{(\nu'_z,\nu'^\vee_{w})} & \A'\times_{S'} \A'^\vee\ar{d}\\
			S\ar{r}{(\nu_z,\nu^\vee_{w})} & \A\times_S \A^\vee
		\end{tikzcd}
	\]
	Let ${\PP'_h}$ resp.\ ${\PP_h}$ be the Poincaré bundles associated to $\A'$ resp.\ $\A$; denote their unique extensions from Theorem \ref{Extension of Poincare bundle} by $\overline{\PP'_h}$ resp.\ $\overline{\PP_h}$; denote the continuously extended canonical metrics by $||.||'$ resp.\ $||.||$ (Theorem \ref{continuous extension of the Poincaré bundle}). To sum up, we have the relation
	\begin{equation}\label{base change of cycles vs normal function}
			\BB_{z',w'} \cong (\nu'_z,\nu'^\vee_{w})^*\PP'_h.
	\end{equation}

	\emph{Claim}. The Lear extension of $\BB_{z,w}$ is given by 
	\[
	[\BB_{z,w},||.||]_{\Sbar}\cong \frac{1}{[k':k]}\nu_*[\BB_{z',w'},||.||']_{\Sbar'}.
	\] 
	Indeed, we observe
	\begin{align*}
		(\nu_*[\BB_{z',w'},||.||']_{\Sbar'})|_S &\cong (\nu|_S)_*([\BB_{z',w'},||.||']_{\Sbar'}|_{S'}) & \text{by (*) below}\\ 
		&\cong (\nu|_S)_*\BB_{z',w'}\\ 
		&\cong (\nu|_S)_*\BB_{<\nu'_z,\nu'^\vee_w>}\\
		&\cong [k':k]\BB_{<\nu_z,\nu^\vee_w>}\\
		&\cong [k':k]\BB_{z,w}.
	\end{align*}
	Here (*) is due to the base change theorem (\href{https://stacks.math.columbia.edu/tag/02RG}{Stacks Project, Tag 02RG}); the second isomorphism is just by definition of Lear extension, and the rest follows from \eqref{Base change and height, formula} and  Corollary \ref{Lear extension and Neron model}.  The line bundle $\nu_*[\BB_{z',w'},||.||']_{\Sbar'}$ admits a continuous metric  induced by the finite morphism $\nu$ and it extends $[k':k]\BB_{z,w}$,  so the claim follows from the uniqueness of the Lear extension. 
	
	Finally, by \eqref{desired commutativity for comparison of Poincare bundles} (plus the argument of  \ref{NT pairing and dual}), we obtain an isomorphism 
	\begin{equation}\label{iso for pullback of P_a and P_h, 1}
		\left(\widetilde{\nu'_{z}},k^p_X\widetilde{\nu'^\vee_{w}}\right)^*\overline{\PP'_h}\cong \left(\widetilde{\theta^p_{z'}},\widetilde{\lambda^q_X\circ \theta^q_{w'}}\right)^*\overline{\PP'_a} \in \Pic(\overline{S}')_\Q,
	\end{equation}
	where $\PP'_a$ is the Poincaré bundles of $\BB\times_S S'$, and the overline indicates the extension provided by Theorem \ref{Extension of Poincare bundle}. We now push forward the isomorphism \eqref{iso for pullback of P_a and P_h, 1} along $\nu\colon \Sbar'\to \Sbar$ (with Corollary \ref{Lear extension and Neron model} in mind) and find 
	\[
	k^p_X[k':k][\BB_{z,w},||.||]_{\Sbar}\cong \nu_*\BB_{<\theta^p_{z'},\lambda^q\circ \theta^q_{w'}>}\cong [k':k]\BB_{<\theta^p_{z},\lambda^q\circ \theta^q_{w}>},
	\] where the first isomorphism is by the \emph{Claim}, the second one is by formula \eqref{Base change and height, formula}. We finish by dividing both sides by $[k':k]$.
\end{proof}
Let us summarize our finding:
\begin{thm}\label{Main Result, alg. trivial cycles}
	Let $X$ be a nice $d$-dimensional over a complex function field $\C(S)$ of a smooth $\C$-curve $S$.  Let $z\in\CH^p_\alg(X), w\in \CH^q_\alg(X)$ with $p+q=d+1$. Let $\X\to \Sbar$ be a regular model of $X\to \Spec k$ over a smooth partial compactification $\Sbar\supset S$. The following  $\Q$-line bundles over $\Sbar$ associated to $z, w$ are isomorphic:
	\begin{enumerate}
		\item[(a)]  the geometric  height pairing bundle $\E_{<z,w>}$; 
		\item [(b)]  the Néron--Tate bundle $(k^p_X)^{-1}\BB_{<\theta^p(z),\lambda^q_X\circ \theta^q(w)>}$.
	\end{enumerate}
Moreover, if Griffiths's Conjecture \ref{Griffiths conjecutre for incidentally trivial cycles} holds for $\CH^p_\alg(X)$, then the two $\Q$-line bundles are isomorphic to 
\begin{enumerate}
	\item[(c)] the Lear extension $[\BB_{z,w},||.||]_{\overline{S}}$   (recall the convention in \ref{spread out normal function}).
\end{enumerate}
\end{thm}
\begin{proof}
	This is a combination of Theorems \ref{Kue96, Thm.8.2} and  \ref{comparison biextension vs NT-height line bundles}.
\end{proof}
\section{Comparison of $\mathbb{G}_m$-biextensions}\label{Section: G_m-biextensions}
The goal is to complete the proof of the commutativity of  Diagram \eqref{normal functions and theta}. Resume the notation there. In fact, we want to prove the following (recall  \ref{K-points of Jac and Pic})
\begin{prop}\label{commutativity of PD_a and h^p-dual}
Assume	Griffiths's conjecture that $h^p$ is an isogeny (Lemma \ref{hp isog. and Griffiths conj}). Then for $w\in CH^q_\alg(X_K)$,  we have
	\[
	h^{p,\vee}\circ \lambda^q_X\circ \theta^q(w)=k^p_X\PD_a\circ \AJ^q(w)\in J^{p,\vee}_{a,X_K}(K).
\]
\end{prop}
Keep assuming Griffiths's conjecture.
\begin{cor}
	$\PD_a\colon J^q_{a,X_K}\to J^{p,\vee}_{a,X_K}$ is an isogeny.
\end{cor}
\begin{cor}
	The homomorphism 
	\[
	k^p_X\PD_a\circ \AJ^q\colon \CH^q_\alg(X_K)\to J^{p,\vee}_{a,X}(K)
	\] is a Picard homomorphism.
\end{cor}
	By Chow's rigidity theorem (Lemma \ref{Chow's rigidity}), it suffices to prove  Proposition \ref{commutativity of PD_a and h^p-dual} for nice $\C$-varieties. Thus from now on assume that $X$ is a nice $\C$-variety and the goal is to prove
\begin{prop}\label{commutativity of PD_a and h^p-dual, complex}
	Assume Griffiths's conjecture. We have 
	\[
		h^{p,\vee}\circ \lambda^q_X\circ \theta^q(w)=k^p_X\PD_a\circ \AJ^q(w)\in J^{p,\vee}_{a,X}(\C),\quad w\in \CH^q_\alg(X).\]
	\end{prop} (As before,  $p+q=d+1$ and $d=\dim X$.)

A clean, conceptual (if somewhat complicated-looking) way to prove the identity is through biextensions.

The notion of $\mathbb{G}_m$-biextensions is introduced by Mumford (\cite{MumfordBiext}) and vastly generalized by Grothendieck in SGA 7. In fact, it can be defined in any topos, and the multiplicative group $\bG_m$ can be replaced by  any abelian group object  $\mathcal{G}$  in that topos. Roughly speaking,  a $\mathcal{G}$-biextension over a product $\calA\times \calB$ of abelian group objects is a $\mathcal{G}$-torsor with two compatible partial group laws.  We will not present the precise definition but refer to the useful references  \cite{BlochBiext}, \cite{MS95}, \cite[\S 10]{BP19} and  the theses \cite{Meyer}, \cite{Seibold}.

We will be using $\C^\times$-biextensions  of abelian groups (in the topos  of sheaves over a one-point set) or $\bG_{m}$-biextension  of  compact complex tori  (in the topos  of sheaves over the big complex-analytic site).

Let $\mathcal{A}, \mathcal{B}$ be abelian group objects in either of the above two topoi. Let $\mathrm{Biext}(\mathcal{A}\times\mathcal{B};\mathbb{G}_m)$ ($\bG_{m,\C}=\C^\times$) denote the set of isomorphism classes of $\mathbb{G}_m$-{biextensions} over $\mathcal{A}\times\mathcal{B}$. The set has a natural abelian group structure.

\begin{void}
If $\calA, \calB$ are (represented by) compact complex tori, then we have an identification
\begin{equation}\label{identification of biextensions}
\mathrm{Hom}(\calA,\calB^\vee)=\mathrm{Biext}(\calB\times\calA;\mathbb{G}_m),\quad f\mapsto (\id_B\times f)^*\PP_B\eqcolon \PP_f
\end{equation} 
where $\PP_B$ is the Poincaré (biextension) bundle on $\calB\times\calB^\vee$. This is the prototype of biextensions.
\end{void}
\begin{void}[Slices]\label{slices of biextensions}
	Let  $P\in \mathrm{Biext}(A\times B;\C^\times)$ be a biextension of abelian groups. For each $b\in B$, we can consider the $b$-\textbf{slice}  $P|_{A\times b}$  by restriction, which yields an exact sequence 
	\[
	0\to \C^\times \to P|_{A\times b}\to A\to 0
	\] of abelian groups. 
	
	Now suppose there is a homomorphism $\phi\colon B'\to B$ of abelian groups. It induces a homomorphism
	\begin{equation}
		(1\times \phi)^*\colon \mathrm{Biext}(A\times B;\C^\times)\to \mathrm{Biext}(A\times B';\C^\times)
	\end{equation}   such that for $b'\in B'$, we have
\begin{equation}
	(1\times \phi)^*P|_{A\times b'}=P|_{A\times \phi(b')}.
\end{equation}
Certainly, a similar construction can be performed in the other factor. A general pullback of biextensions along a  homomorphism $\psi\times \phi\colon A'\times B'\to A\times B$ may be understood ``slice-wise'' by writing $\psi\times \phi=(\psi\times 1) \circ (1\times \phi)$ (pullback has the expected functoriality).
\end{void}

\begin{void}
	Let $\PP_{\PD}\in \mathrm{Biext}(J^p_X\times J^q_X;\bG_m)$ be the Poincaré biextension corresponding to the Poincaré duality $\PD\colon J^q\cong J^{p,\vee}$ (notation as in \eqref{identification of biextensions}). It can be verified that  $\PP_{\PD}$ (in the complex-analytic category) induces a $\C^\times$-biextension 
	\[
	\P^\hom\to J^p_{X}(\C)\times J^q_X(\C)
	\]  of abelian groups (see e.g.\ \cite[9.4]{Seibold}).

On the other hand, Bloch constructed a $\C^\times$-biextension $$\bE^\hom\to \CH^p_\hom(X)\times \CH^q_\hom(X)$$ and conjectured (at the end of \cite{BlochBiext}) that there is a natural isomorphism of $\C^\times$-biextensions over $J^p(\C)\times J^q(\C)$:
	\begin{displaymath}
		(\AJ^p\times \AJ^q)^*\P^\hom\cong \bE^\hom,
	\end{displaymath} which is proved by Müller-Stach (\cite[Theorem 1]{MS95}). The slices of the two biextensions are related by the \emph{pullback} diagram (in the category of abelian groups):
\begin{equation}
	\begin{tikzcd}
		0\ar{r} & \C^\times \ar[equal]{d}  \ar{r}& \bE^\hom|_{\CH^p_\hom(X)\times w} \ar{d} \ar{r} & \CH^p_\hom(X) \ar{d}{\AJ^p}\ar{r} & 0\\
		0\ar{r} & \C^\times \ar{r} & \P^\hom|_{J^p_X(\C)\times \AJ^q(w)} \ar{r}& J^p_X(\C) \ar{r} & 0
	\end{tikzcd}   
\end{equation} for $w\in \CH^q_\hom(X)$.
See also O. Meyer's thesis \cite{Meyer} for a detailed exposition.
\end{void}
\begin{void}
Let 
\[
\bE\to \CH^p_\alg(X)\times \CH^q_\alg(X)
\] be the $\C^\times$-biextension restricted from $\bE^\hom$. Let 
\[
\P\to J^p_{a,X}(\C)\times J^q_{a,X}(\C)
\]
be the $\C^\times$-biextension induced by $\PP_{\PD_a}$.  Then Müller-Stach's theorem implies
\begin{equation}\label{MüllerStach, algebraic}
	(\AJ^p\times\AJ^q)^*\P\cong \bE.
\end{equation} 
\end{void}

\begin{void}
	One of the main results of M. Seibold’s (PhD) thesis \cite{Seibold} is to relate the right hand side of \eqref{MüllerStach, algebraic} to the Poincaré biextension bundles of higher Picard varieties $\Pic^\bullet_X$, which is of algebraic nature.
	
	Let $\PP_{\lambda^q}\in \mathrm{Biext}(\Pic^p_X\times \Pic^q_X;\bG_m)$, where $\lambda^q\coloneq\lambda^q_X\colon \Pic^q_X\to \Pic^{p,\vee}_X$ is the isogeny from \eqref{Pic(alpha^*)}. We obtain an induced $\C^\times$-biextension
	\[
	\P_\lambda\to \Pic^p_X(\C)\times \Pic^q_X(\C).
	\]
\end{void}
\begin{void}
Rearranging the defining property of Poincaré-$p$-cycle $\alpha^p\coloneq \alpha^p_X\in \CH^p(\Pic^p_X\times X)$ (\ref{Pic of Poincare p-cycle}), we see that there is a commutative diagram
	 \begin{equation}
	 	\begin{tikzcd}\label{theta-p-q}
	 		\CH^\alg_0(\Pic^p_X)\times \CH^q_\alg(X) \ar{r}{\theta_0\times \theta^q} \ar[swap]{d}{\alpha^p_*\times 1}& \Pic^p_X(\C)\times \Pic^q_X(\C)\ar{d}{[k^p_X]\times 1}\\
	 		\CH^p_\alg(X)\times \CH^q_\alg(X) \ar{r}{\theta^p\times \theta^q}& \Pic^p_X(\C)\times \Pic^q_X(\C)
	 	\end{tikzcd}
	 \end{equation} where the $\theta_0$ in the above diagram arises from the (universal) Albanese map
 \[
 \theta_0=\theta^{e}\colon \CH^{e}_\alg(\Pic^p_X)\to \Pic^p_X(\C)
 \] under the canonical identification $\Pic^e_{\Pic^p_X}=\Pic^p_X$; here  $e\coloneq\dim \Pic^p_X$.
\end{void}
\begin{prop}\label{Seibold, Satz 10.12}
	There is an isomorphism $\C^\times$-biextensions
	\begin{equation}
		(\alpha^p_*\times 1)^*\bE=(\theta_0\times \theta^q)^*\P_{\lambda}\in \mathrm{Biext}(\CH^\alg_0(\Pic^p_X)\times \CH^q_\alg(X);\C^\times).
	\end{equation}
\end{prop}
\begin{proof}
	See \cite[Satz 10.12]{Seibold}.
\end{proof}
To ease notation, and also because a morphism of complex abelian varieties is determined by its induced homomorphism on $\C$-points, we will not distinguish between a  morphism and the induced group homomorphism.
\begin{void}
	From \eqref{MüllerStach, algebraic}, \eqref{theta-p-q}, Proposition \ref{Seibold, Satz 10.12} and the compatibility $\theta^\bullet=h^\bullet\circ \AJ^\bullet$, we conclude that
	\begin{equation}\label{intermediate conclusion of biextensions}
		(\alpha^p_*\times 1)^*(\AJ^p\times \AJ^q)^*(h^p\times h^q)^*\P_\lambda\cong (\alpha^p_*\times 1)^*(\AJ^p\times\AJ^q)^*([k^p_X]\times 1)^*\P,
	\end{equation}
	as $\C^\times$-biextensions over $\CH^\alg_0(\Pic^p_X)\times \CH^q_\alg(X)$. Note that we can freely interchange the order of the pullback along $[k^p_X]\times 1$. 
\end{void}
\begin{void}[Goal]\label{iso of biextension, final}
	We want to prove that
	\[
	(h^p\times h^q)^*\PP_{\lambda^q}\cong ([k^p_X]\times 1)^*\PP_{\PD_a}
	\] as $\bG_m$-biextensions over $J^p_{a,X}\times J^q_{a,X}$. This will imply Proposition \ref{commutativity of PD_a and h^p-dual, complex} (under the identification \eqref{identification of biextensions}). 
\end{void}
\begin{void}
	Observe that  the composition
	\[
	\CH^\alg_0(\Pic^p_X)\xrightarrow{\alpha^p_*} \CH^p_\alg(X)\xrightarrow{\AJ^p} J^p_{a,X}(\C)
	\] is a regular homomorphism. Since $\Pic^p_X=\mathrm{Alb}_{\Pic^p_X}$ is the algebraic representative on $\CH^\alg_0(\Pic^p_X)$, there is a unique morphism $f^p\colon \Pic^p_X\to J^p_{a,X}$ such that
	\begin{equation}\label{factorization through alb}
		\AJ^p\circ \alpha^p_*=f^p\circ\theta_0.
	\end{equation}
Moreover, by \eqref{theta-p-q} and \eqref{factorization through alb}, we have
\[
h^p\circ f^p\circ \theta_0=h^p\circ \AJ^p\circ \alpha_{*}^p=\theta^p\circ \alpha^p_*=[k^p_X]\circ \theta_0.
\] Since $\theta_0$ is surjective, we find
\begin{equation}\label{hp-fp is kp}
	h^p\circ f^p=[k^p_X]\in \mathrm{End}(\Pic^p_X).
\end{equation}
\end{void}
\begin{lem}
The homomorphism 
	\[
	(\theta_0\times 1)^*\colon \mathrm{Biext}(\Pic^p_X(\C)\times \Pic^q_{a,X}(\C);\C^\times)\to \mathrm{Biext}(\CH_0^\alg(\Pic^p_X)\times \Pic^q_{a,X}(\C);\C^\times)
	\] is injective.
\end{lem}
\begin{proof}
(The proof works for any complex abelian variety than just for $\Pic^\bullet_X$.)	Let $0\in \Pic^p_X(\C)$ denote the identity element. The map $\theta_0$ has a section 
	\[
	  \Pic^p_X(\C)\to \CH_0^\alg(\Pic^p_X),\quad x\mapsto [x]-[0].
	\]  This implies the injectivity of $(\theta_0\times 1)^*$.
\end{proof}
\begin{void}
	Consider the homomorphism
	\begin{displaymath}\label{final injectivity of pullback of biextensions}
		(f^p\times \AJ^q)^*\colon \mathrm{Biext}(J^p_{a,X}(\C)\times J^q_{a,X}(\C);\C^\times)\to \mathrm{Biext}(\Pic^p_X(\C)\times \CH^q_\alg(X);\C^\times).
	\end{displaymath}
Using \eqref{identification of biextensions} and \ref{slices of biextensions}, we can interpret the above pullback as a composition of homomorphisms
	\begin{equation}\label{injectivity check, final}
		\CH^q_\alg(X)\xrightarrow{\AJ^q} J^q_{a,X}(\C)\xrightarrow{\varphi(\C)}J^{p,\vee}_{a,X}(\C)\xrightarrow{(f^p)^\vee(\C)}\Pic^{p,\vee}(\C),
	\end{equation}
	 where $\varphi$ is the morphism $\varphi\colon J^q_{a,X}\to J^{p,\vee}_{a,X}$ corresponding to the difference biextension
	\[
	k^p_X\PP_{\PD_a}-(h^p\times h^q)^*\PP_{\lambda^q}\in \mathrm{Biext}(J^p_{a,X}\times J^q_{a,X};\bG_m),
	\]
By assumption, the composition \eqref{injectivity check, final} is the zero map. We want to conclude  $\varphi(\C)\equiv 0$, which implies the triviality of the difference  biextension.
	
Now recall \eqref{hp-fp is kp}. If Griffiths conjecture is true, or equivalently (Lemma \ref{hp isog. and Griffiths conj}), if $h^p$ is an isogeny, then $f^p$ is also an isogeny (since $[k^p_X]$ is), in particular $f^p$ is surjective. Thus $(f^p)^\vee\colon J^{p,\vee}_{a,X}\to \Pic^{p,\vee}_X$ is injective, hence it is also injective on $\C$-points. 

By assumption, we have
\[
f^{p,\vee}(\C)\circ \varphi(\C)\circ \AJ^q(w)=0
\] for all $w\in \CH^q_\alg(X)$.  Since $\AJ^q$ is surjective, we obtain $f^{p,\vee}(\C)\circ \varphi(\C)\equiv 0$. By the injectivity of $f^{p,\vee}(\C)$, it follows that $\varphi(\C)\equiv 0$.
\end{void}
\bibliographystyle{alpha}
\bibliography{RefforPhDStudy}
\end{document}